\documentclass[12pt,oneside]{amsart}
\usepackage{graphicx}
\usepackage{amssymb}
\usepackage{epstopdf}
\usepackage{typearea}
\usepackage{amscd,amssymb,verbatim,array}
\usepackage{hyperref}
\usepackage{charter}
\usepackage{pdfsync}
\usepackage{color}
\usepackage{mathrsfs}
\numberwithin{equation}{section} 
\theoremstyle{plain}
\usepackage[a4paper,left=2.5cm,right=2.5cm, top=2.5cm,bottom=2cm]{geometry}

\newcommand{\field}[1]{\mathbb{#1}}

\newcommand{\R}{\field{R}}
\newcommand{\C}{\field{C}}
\newcommand{\N}{\field{N}}

\newcommand{\E}{\field{E}}
\renewcommand{\P}{\field{P}}

\newcommand{\cali}[1]{\mathscr{#1}}

\newcommand{\B}{\mathbb{B}}
\newcommand{\D}{\mathbb{D}}
\renewcommand\P{\mathbb{P}}
\newcommand{\Leb}{\mathop{\mathrm{Leb}}\nolimits}
\renewcommand{\S}{\mathbb{S}}

\DeclareMathOperator\supp{supp}
\DeclareMathOperator\Div{Div}
\DeclareMathOperator\ddc{dd^c}
\DeclareMathOperator\dist{dist}

\newtheorem{theorem}{Theorem}[section]
\newtheorem{definition}[theorem]{Definition}
\newtheorem{proposition}[theorem]{Proposition}

\newtheorem{lemma}[theorem]{Lemma}
\theoremstyle{definition}
\newtheorem{remark}[theorem]{Remark}

\newtheorem{example}[theorem]{Example}

\newtheorem{problem}[theorem]{Problem}

\title[Asymptotic equilibrium distribution of zeros of random holomorphic sections]
{A survey on asymptotic equilibrium 
distribution of zeros of random holomorphic sections}

\author[George Marinescu]{George Marinescu}
\address{Universit{\"a}t zu K{\"o}ln,  Mathematisches Institut,
		Weyertal 86-90, 50931 K{\"o}ln, Germany
		\newline\mbox{\quad}\,Institute of Mathematics `Simion Stoilow', 
		Romanian Academy, Bucharest, Romania}
\thanks{Partially supported by the DFG funded projects
SFB TRR 191 `Symplectic Structures in Geometry, 
Algebra and Dynamics' (Project-ID 281071066\,--\,TRR 191),
DFG Priority Program 2265 `Random Geometric Systems' (Project-ID 422743078), the ANR-DFG project `QuasiDy\,--\,Quantization, Singularities, 
and Holomorphic Dynamics' (Project-ID 490843120) and the DFG Project `Intersection of closed positive currents´  (Projektnummer-500055552)}

\author[Duc-Viet Vu]{Duc-Viet Vu}
\address{Universit{\"a}t zu K{\"o}ln,  Mathematisches Institut,
		Weyertal 86-90, 50931 K{\"o}ln, Germany}

\date{March 3, 2025} 

\begin{document}
\begin{abstract}
This is a survey of results concerning the asymptotic equilibrium 
distribution of zeros of random holomorphic polynomials and holomorphic
sections of high powers of a positive line bundle, as related to \cite{Marinescu-Vu}.
Our primary focus is on the role of pluripotential theory in this research area. 
\end{abstract}

\maketitle
\tableofcontents

\section{Introduction}

The objective of this survey is to review results
pertaining to the asymptotic equilibrium 
distribution of zeros of random holomorphic sections of high
powers of a positive line bundle, 
as related to the authors' paper on
Bergman kernel functions associated to
measures supported on totally real submanifolds \cite{Marinescu-Vu}.
These results encompass generalizations of the notion of equilibrium potential and measure
and are closely related to
asymptotics of Bergman kernels relative to a wide range of weights and measures.
%We will be mainly be concerned with the area around the paper \cite{Marinescu-Vu}.

A basic example is that of Kac polynomials, that is, polynomials 
$$f(z)=\sum_{j=0}^k a_{kj}z^j$$ with 
independent and identically distributed (i.i.d.) standard Gaussian coefficients 
$a_{kj}$. 
The asymptotic equilibrium distribution of zeros for Kac polynomials is
described by the classical result of Kac \cite{Kac43} and 
Hammersley \cite{Ham56}. This result states that as the degree of 
the polynomial $k$ approaches infinity, the zeros of the polynomial 
tend to accumulate on the unit circle $\S^1\subset\C$. 
For this purpose, we identify the space $\mathcal{P}_k$ of polynomials in 
$z$ with degree at most $k$ with the space $\C^{k+1}$ by means 
of the basis $\{1, z, \ldots, z^k\}$. We then endow $\mathcal{P}_k$ 
with a Gaussian probability measure originating from $\C^{k+1}$.
A conceptual and very fertile manner to prove the above result 
\cite{Bloom05,Bloom09,BlLe15,ShZ03} is to
introduce the following inner product on  the space $\mathcal{P}_k$ 
$$(f,g)=\frac{1}{2\pi}\int_{\S^1}f(e^{i\theta})
\overline{g(e^{i\theta})}d\theta,\quad f,g\in\mathcal{P}_k,$$ 
for which an orthonormal basis is $1, z,\ldots,z^k$.
The key of the proof is the fact that the Bergman kernel function
$B_k(z)$ of $(\mathcal{P}_k,(\cdot,\cdot))$ satisfies 
$$\frac{1}{2k}\log B_k(z)\to\log^{+}|z|,\quad
\text{locally uniformly on $\C$ as $p\to\infty$}$$
and $\frac{1}{2\pi}\log^{+}|z|$ is a potential of the delta function 
$\delta_{\S^1}$, which is also the equilibrium measure of the unit circle $\S^1$.

%Bloom \cite{Bloom05} highlighted the importance of extremal subharmonic functions
%in the asymptotic equilibrium distribution result.
%He also introduced Bernstein-Markov measures as the most comprehensive 
%framework for defining Gaussian random sections.

Another basic example is that of $SU(2)$ polynomials. 
They serve as a basic model for eigenfunctions 
of quantum chaotic Hamiltonians and asymptotic distribution of their zeros has 
been intensively studied by physicists. 
The zeros of Gaussian $SU(2)$ polynomials 
\[f(z)=\sum_{j=0}^k a_{kj}\,\sqrt{\binom{k}{j}}\,z^j\]
become uniformly distributed over the Riemann sphere with respect to the 
Fubini-Study volume as the degree $k$ tends to infinity. 

In the same vein, one can view zeros of polynomials in one complex variable as 
interacting particles, as, for instance, eigenvalues of random asymmetric matrices 
can be physically interpreted as a two-dimensional electron gas confined in a disk. 
In this manner, a closely related model is Weyl polynomials 
\[f(z)=\sum_{j=0}^k a_{kj}\,\frac{1}{\sqrt{j!}}\,z^j\] 
whose zeros become uniformly distributed with respect to the Lebesgue 
measure in the unit disc.

Let us now describe the general setting.
Let $K$ be a non-pluripolar set in $\C^n$ and let $\mu$ 
be a probability measure on $\C^n$ such that the support of $\mu$
is contained in $K$ and is non-pluripolar. 
Let $Q$ be a continuous weight on $K$. 
We will assume that the triple $(K,Q,\mu)$ is Bernstein-Markov, 
cf.\ Definition \ref{D:BeMa}. 
Let $\mathcal{P}_k(K)$ 
be the space of restrictions of complex polynomials of degree at most $k$ 
in $\C^n$ to $K$. Let $d_k:= \dim \mathcal{P}_k(K)$, 
and let  $s_1,\ldots, s_{d_k}$ be an orthonormal basis of 
$\mathcal{P}_k(K)$ with respect to the $L^2(\mu,k Q)$-scalar product
\eqref{eq:l2muk}.  
Consider the random polynomial
\begin{align}\label{eq-randompoly}
p_k:= \sum_{j=1}^{d_k} a_{kj} s_j,
\end{align}
where $a_{kj}$ are complex i.i.d. random variables. 
%The study of zeros of random polynomials has a long history. 
A classic example is the Kac polynomial, where $n=1$ and $p_j=z^j$. 
A more general setting involves replacing $\mathbb{C}^n$ with 
a compact projective manifold $X$ and $K$ as a non-pluripolar set in $X$. 
Polynomials of degree $k$ are replaced by sections of $L^k$, 
where $L$ is a positive line bundle on $X$.
In this more general setting, there has been a substantial amount of work, 
primarily focusing on the case where $K=X$. In that setting, under mild assumptions on the distribution of random coefficients $a_{kj}$, it was shown that the (random) normalized zero divisor of $p_k$ converges weakly to an equilibrium current almost surely. This property is called the equidistribution of zeros of random polynomials toward the equilirium current.  
For a concise overview of this aspect of the problem, we refer to the survey 
\cite{Bayraktar-Coman-Hendrik-M}, and also to 
\cite{Drewitz-M-L3,Drewitz-M-L2,Drewitz-M-L1} for recent work in 
the case where $X$ is non-compact. One can furthermore consider 
holomorphic sections with prescribed vanishing order along a divisor, 
in which case, there is a corresponding notion for the equilibrium current. 
We refer to \cite{Comna-M-N-vanishinalong3,Comna-M-N-vanishinalong2,Comna-M-N-vanishinalong} 
for details.

As we shall see, quantitative equidistribution
has not been studied extensively, except for the case where $K=X$. 
In this case, quantitative estimates have been obtained due to the smooth 
asymptotic behavior of the Bergman kernel function for positive holomorphic line bundles
\cite{Marinescu-Ma}. In addition, \cite[Corollary 7.4]{DS_tm} gives a quantitative 
result for real polynomials. 
More recently, the quantification of the equidistribution of zeros of random 
polynomials for a relatively large class of compact sets $K$ has been achieved 
in \cite{Marinescu-Vu}. This progress is due to advances in the asymptotic 
behavior of Bergman kernel functions, as described in 
\cite{Berman-OrtegaCerda,Marinescu-Vu}.  

The distribution of zeros of more general random polynomials associated to
orthonormal polynomials (as in \eqref{eq-randompoly}) 
was considered by Shiffman-Zelditch \cite{Shiffman-Zelditch} 
by observing that  $1, z, \ldots, z^k$ form an orthonormal basis of 
the restriction  of the space of polynomials in $\C$ to $\S^1$ with respect to the $L^2$-norm induced 
by the Haar measure $\mu_0$ on $\S^1$.   In this setting, the necessary and 
sufficient conditions for the distribution of $a_{kj}$ so that the zeros of 
$p$ is equidistributed almost surely or in probability with respect to the equilibrium 
measure associated to $(K,Q,\mu)$ as $k \to \infty$ are known if $Q=0$; 
see \cite{Bloom-Dauvergne,Dauvergne,Hammersley,Ibragimov-Zaporozhets}. 
It was also proved in \cite{Shiffman-Zelditch} that the correlations between 
zeros have a universal scaling limit (when normalized properly) 
if $K$ is a simply connected domain with analytic boundary or the 
boundary of such domains. 

We focus in this survey only on results in the higher-dimensional setting.
In Section \ref{S:Bkf} we review the notion of
Bergman kernel function and introduce the Bernstein-Markov measures.
In Section \ref{S:Epsh}, we discuss  extremal plurisubharmonic functions and their regularity.  This notion plays a crucial role in the study of Bergman kernel functions associated with measures supported on real submanifolds.
In Section \ref{S:Zrp} we present a quantiative large deviation type estimate for the equidistribution of zeros of polynomials from \cite{Marinescu-Vu} and propose some open problems in higher dimensions. We also prove Theorems \ref{the-expectedzeros} and \ref{the-zeros-variant} which are variants of results from \cite{Marinescu-Vu}.

\section{Bergman kernel functions}\label{S:Bkf}

In the theory of equidistributon of zeros of holomorphic sections 
of lines bundles, the notion of Bergman (Christoffel-Darboux) kernel function 
plays an essential role.

The Christoffel-Darboux kernel in complex dimension one is a 
classical object of great importance in approximation theory and random matrix theory. 
An extensive literature exists on the asymptotic behavior of the 
Christoffel-Darboux kernel function, whose study can be traced back to 
the early twentieth century \cite{Szegobook}.
We refer to \cite{Beckermann-Putinar-Saff-Sty,BoucksomBermanWitt,
Bloom-Levenberg-PW,Danka-Totik,Deift-book,Dunkl-Xu-book,Kroo-Lubinsky,
Lubinsky,Totik,Xu-simplex} for some recent works on this topic.
The Christoffel-Darboux kernel can be regarded as a specific case 
of the Bergman kernel function in complex geometry, as discussed 
in \cite{BoucksomBerman,Marinescu-Ma}.
The study of the Christoffel-Darboux kernel in higher dimensional 
settings has to be better understood and remains an area of ongoing research.

In this survey, we focus on elucidating the role of the Bergman 
(Christoffel-Darboux) kernel function in the problem of equidistribution of zeros 
of random holomorphic sections. We start with a classical setting in 
$\mathbb{C}^n$ and then explore a more general situation in complex geometry. 
We assume that the reader is familiar with the basic concepts of 
pluripotential theory and complex geometry as outlined in 
the references \cite{Demailly_ag,Klimek}.

Let $K$ be a non-pluripolar compact subset in $\C^n$, 
\emph{i.e,} $K$ is not contained in $\{\varphi = -\infty\}$ 
for any plurisubharmonic (psh) function $\varphi$ on $\C^n$,
which is not identically $-\infty$. 
Let $\mu$ be a probability measure whose support is 
non-pluripolar and is contained in $K$, 
and $Q$ be a real continuous function on $K$. 
Let $\mathcal{P}_k$ be the space of restrictions to $K$ 
of complex polynomials of degree at most $k$ on $\C^n$. 
The  scalar product
\begin{equation}\label{eq:l2muk}
\langle s_1, s_2\rangle_{L^2(\mu,k Q)}:=
\int_K s_1\overline s_2  e^{-2 k Q} d \mu
\end{equation}
induces the $L^2(\mu, k Q)$-norm on $\mathcal{P}_k$.
\emph{The Bergman kernel function of order $k$} 
associated with $\mu$ with weight $Q$ is defined by 
\begin{equation}\label{eq:Bkf1}
B_k:K\to [0,\infty),\quad B_k(x):= \sup_{s \in \mathcal{P}_k} 
|s(x) e^{-k Q(x)}|^2/\|s\|^2_{L^2(\mu, k Q)}, \quad x\in K.
\end{equation} 
Equivalently if $(s_1, \ldots, s_{d_k})$ (here  $d_k$ denotes 
the dimension of $\mathcal{P}_k$) is an orthonormal basis of 
$\mathcal{P}_k$ with respect to the $L^2(\mu,k Q)$-norm given
by \eqref{eq:l2muk}, then 
\begin{equation}\label{eq:Bkf2}
B_k(x)= \sum_{j=1}^{d_k} |s_j(x)|^2 e^{-2k Q(x)},\quad x\in K.
\end{equation}
When $Q \equiv 0$, we say that $B_k$ is \emph{unweighted}. 
In this case the inverse of $B_k$ is known as 
\emph{the Christoffel function} in the literature on orthogonal polynomials.  
In practice we also use a modified version of the Bergman kernel function as follows:
\begin{equation}\label{eq:Bkf3}
\widetilde{B}_k:\C^n\to[0,\infty),\quad
\widetilde{B}_k(x):= \sup_{s \in \mathcal{P}_k} |s(x)|^2/ \|s\|^2_{L^2(\mu, k Q)},
\quad x \in \C^n.
\end{equation} 
The advantage of using $\widetilde{B}_k$ instead of ${B}_k$
is that $\widetilde{B}_k$ is well defined on $\C^n$.  
  
The asymptotics of the Bergman kernel function 
(or its inverse, the Christoffel function) is essential for many applications including approximation theory, 
random matrix theory, etc. There is an immense literature on such asymptotics. 
We refer to \cite{Beckermann-Putinar-Saff-Sty,BoucksomBermanWitt,
Bloom-Levenberg-PW,Danka-Totik,Deift-book,Dunkl-Xu-book,Kroo-Lubinsky,
Lubinsky,Totik,Xu-simplex}, to cite  just few,
for an overview on this very active research field.
%to help readers to capture a rough overview on this very active research field.       

In applications, it is important to consider a large class of measures $\mu$.
Standard examples are measures supported on domains 
on $\R^n\subset\C^n$,  \emph{i.e.}, balls or simplexes in $\R^n$  or in the unit ball in $\C^n$. In dimension one, Christoffel-Darnoux kernels for measures with support on Jordan curves in $\C$ were discussed in  
\cite{Beckermann-Putinar-Saff-Sty,Gustafsson-Putinar-Saff-Stylianopoulos,
Totik-transaction,Totik}. This is a particular case of the case where the measure $\mu$ is supported on piecewise-smooth domains in a generic Cauchy-Riemann submanifold $K$ in $\C^n$. The notion of Bergman (Christoffel-Darboux) kernel function can also be considered in the more general setting of projective manifolds as follows.

%Present basics things about line bundles, positive Hermitian metrics, enveloped, basic properties.

Let $X$ be a projective manifold of dimension $n$.  
Let $(L,h_0)$ be an ample line bundle equipped with 
a Hermitian metric $h_0$ whose Chern form $\omega$ is positive. 
Let $K$ be a compact non-pluripolar subset in $X$.  
Let $\mu$ be a probability measure on $X$ such that the support 
of $\mu$ is non-pluripolar and is contained in $K$. 
Let $h$ be a Hermitian metric on $L|_{K}$ such that 
$h= e^{-2\phi} h_0$, where $\phi$ is a continuous function on $K$. We denote by $H^0(X,L)$ the space of holomorphic sections of $L$. 
For $s_1,s_2 \in H^0(X,L)$, we define
$$\langle s_1, s_2 \rangle:= \int_X \langle s_1, s_2 \rangle_h d \mu.$$
Since $\supp \mu$ is non-pluripolar, the last scalar product defines the $L^2(\mu,h)$-norm on $H^0(X,L)$.  Let $k \in \N$. 
We obtain induced Hermitian metric $h^k$ on $L^k$ and a similar 
$L^2(\mu, h^k)$-norm on $H^0(X,L^k)$.  Put $d_k:= \dim H^0(X,L^k)$. 
Let $\{s_1, \ldots, s_{d_k}\}$ be an orthonormal basis of $H^0(X,L^k)$ 
with respect to $L^2(\mu, h^k)$-norm. 
The \emph{Bergman kernel function of order $k$ associated with $(L,h,\mu)$} is 
\begin{equation}\label{eq:Bkf4}
B_k(x)=B_{k,L,h}(x):= \sum_{j=1}^{d_k} |s_j(x)|_{h^k}^2=
\sup_{s \in H^0(X,L)} |s(x)|_{h^k}^2/ \|s\|_{L^2(\mu,h^k)}
\end{equation}
for $x\in K$.
 
%Compared to works in complex geometry, $e^{-2k \phi} B_k$ is the Bergman kernel function associated 
%to $(L,h_0^k,\mu)$ (which is only defined on $K$ for $\phi$ is only defined on $K$). We stick to 
%$B_k$ rather than $e^{-2 k \phi}B_k$ because it is more convenient for us to have  $B_k$ 
%defined on the whole $X$ in our application later. However we stress that this is purely a notational 
%choice, and is not essential.  

When $\mu$ is a volume form on $X$ and $h=h_0$, 
the Bergman kernel function is the standard Bergman kernel function which is an object of great importance 
in complex geometry, see \cite{Marinescu-Ma} 
for a comprehensive study.

The setting considered previously in $\C^n$ corresponds to the case where 
$X= \P^n$ and $(L,h_0)=(\mathcal{O}(1),h_{FS})$ is the 
hyperplane line bundle on $\P^n$ 
endowed with the Fubini-Study metric. We consider $\C^n$ 
as an open subset in $\P^n$ and the weight $Q$ corresponds to 
$\phi + \frac{1}{2}\log (1+ |z|^2)$. %Recall that there is a natural 
%identification between $\mathcal{L}(\C^n)$ and the set of $\omega_{FS}$-psh functions on $\P^n$ (where $\omega_{FS}$ denotes the Fubini-Study form on $\P^n$) given by 
%\begin{align} \label{corres-Leleongclass}
%u \, \longleftrightarrow \, u- \frac{1}{2}\log(1+ |z|^2)\,,\quad 
%u\in \mathcal{L}(\C^n)\,.
%\end{align}
Another well-known example is the case where $K$ is the unit sphere in 
$\R^n$ (here $n\ge 2$; see, e.g, \cite{Marzo-OrtegaCerda}) 
and $X$ is the complexification of $K$, \emph{i.e,} 
$K= \S^{n-1} \subset \R^n$ which is considered as usual 
a compact subset of $X:= \{z_0^2+ z_1^2+ \cdots+ z_n^2=1\} \subset \P^n$. The line bundle $L$ on 
$X$ is the restriction of the hyperplane bundle $\mathcal{O}(1)\to\P^n$ to $X$. 
We remark that in this case $H^0(X, L^k)$ is equal to the restriction 
of the space of $H^0(\P^n, \mathcal{O}(k))$ to $X$. 
Hence  the restriction of $H^0(X,L^k)$ to $K$ is that of the space of 
complex polynomials in $\C^n$ to $K$.  
%To see this, notice that $X$ is a smooth hypersurface in $\P^n$. 
%Consider the standard exact sequence of sheaves:
%$$0 \to \mathcal{O}(k- \deg X)\to\mathcal{O}(k)\to 
%\mathcal{O}(k)|_X \rightarrow 0,$$
%where the second arrow is the multiplication by a section of 
%$\mathcal{O}(\deg X)$ whose zero divisor is equal to $X$.  
%We thus obtain a long exact sequence of cohomology
%spaces: 
%\begin{multline*}
%0 \to H^0(\P^n,\mathcal{O}(k - \deg X)) \to 
%H^0(\P^n,\mathcal{O}(1)) \to  H^0(\P^n,\mathcal{O}(k)|_X)\\
% \to H^1(\P^n,\mathcal{O}( k- \deg X)) \to \cdots
% \end{multline*}
%In this sequence, $H^0(\P^n,\mathcal{O}(k)|_X)$is isomorphic to $H^0(X,\mathcal{O}(k)|_X)$, and by the Kodaira-Nakano vanishing theorem we have$H^1(\P^n,\mathcal{O}(k- \deg X))= 0$; see \cite[p.\ 156]{Griffiths-Harris}.
As above, the weight $Q$ on $K$ in the spherical model corresponds to 
$\phi+ \frac{1}{2}\log (1+ |z|^2)|_X$ in the setting $(K,X, \mathcal{O}(1)|_X)$.

%In view of asymptotic of Bergman kernel functions for volume forms in $X$, our goal is to obtain bounds of similar flavor for Bergman kernel functions associated to measures supported on totally real (or more generally generic Cauchy-Riemann) submanifolds of $X$. As an application, we obtain a quantitative large deviation type result for zeros of random polynomials. To state our results, we need to recall some notions.    

\begin{definition}\label{D:BeMa} (\cite{BoucksomBerman,Bloom05,Bloom09})
The measure $\mu$ is said to be a \emph{Bernstein-Markov measure} 
(with respect to $(K,\phi,L)$) if  for every constant $\varepsilon>0$ 
there exists a constant $C=C(\varepsilon)>0$ such that  
\begin{align}\label{ine-def-BMintro}
\sup_K |s|_{h^k}^2 \le C e^{\varepsilon k} \|s\|^2_{L^2(\mu, h^k)}
\end{align}
for every $s \in H^0(X, L^k)$. In other words, the Bergman kernel function 
of order $k$ grows at most subexponentially, \emph{i.e,}  
$\sup_K B_k =O(e^{\varepsilon k})$ as $k \to \infty$ for every $\varepsilon>0$.  If  $\mu$ is a Bernstein-Markov measure with respect to $(K,\phi,L)$, we also say that the triple $(K,\mu,\phi)$ (when $L$ is clear from the context) is \emph{Bernstein-Markov or satisfies the Bernstein-Markov property.}

\end{definition}
 
The notion of Bernstein-Markov measure (in the higher dimensional setting) 
was introduced in \cite{Bloom05,Bloom09}.  We refer to \cite{Bloom-Levenberg-PW} 
for a survey about these measures. In the literature on orthogonal polynomials 
(in dimension one), Bernstein-Markov measures are called measures 
in \textbf{Reg} class, see \cite{Stahl-Totik-generalop}. We refer to Theorem \ref{the-BMproperty} below for a large class of Bernstein-Markov measures.

\section{Extremal plurisubharmonic functions}
\label{S:Epsh}

\subsection{Definition and continuity of extremal plurisubharmonic functions}
The second crucial ingredient in the study of the equidistribution of zeros of random holomorphic sections of line bundle is the notion of extremal plurisubharmonic envelopes. We discuss now its properties. It is emphasized that both Bergman kernel function and extremal psh envelope are expected to have other important applications, e.g., in approximation theory. 

%The notion of (Siciak-Zahariuta) extremal plurisubharmonic functions (CITE NGOC CUONG NGUYEN PAPER) plays a central role in the study of equidistribution of zeros of random holomorphic sections as well as related problems. 

We start with the local setting on $\C^n$. Let $K$ be a non-pluripolar compact set in $\C^n$. Recall that for subsets in $\C^n$ the local pluripolarity is the same as the (global) pluripolarity. Let $Q$ be a real continuous function on $K$. We recall that
$$V_{K,Q}:= \sup\{\psi \in \mathcal{L}(\C^n): \, \psi \le Q \quad 
\text{ on } K\},$$
where $\mathcal{L}(\C^n)$ is the set of psh functions $\psi$ on $\C^n$
such that $\psi(z)- \log |z|$ is bounded at infinity on $\C^n$.
If $Q \equiv 0$, we put $V_K:= V_{K,0}$. The function $V^*_{K,Q}$ is a called \emph{the extremal plurisubharmonic function associated with $(K,Q)$}. %It was introduced by ??? 

\begin{lemma} \label{le-lowersemi}
 The upper semi-continuous regularisation 
$V_{K,Q}^*$ of $V_{K,Q}$ belongs to  $\mathcal{L}(\C^n)$. The function $V_{K,Q}$ is always lower semi-continuous. %These two properties also hold if $K$ is only bounded and non-pluripolar.
\end{lemma}

\proof The first assertion follows from the non-pluripolarity of $K$ and the fact that $K$ is compact. The second claim is a consequence of  Lemma \ref{le-pointwiseconvergence} below  or \cite[Corollary 5.1.3]{Klimek}.
\endproof
 We note that the above definition of $V_{K}$  still make sense if $K$ is non-pluripolar bounded subset in $\C^n$.
 %above definition it suffices to have that  $K$ is non-pluripolar, the compactness of $K$ is  

\begin{definition} \label{def-localregular} We say that $(K,Q)$ is \emph{regular} if $V_{K,Q}= V^*_{K,Q}$ (or equivalently, $V_{K,Q}$ is continuous).  The set $K$ is said to be \emph{locally regular} 
if for every $z\in K$,  every open subset $U'$ in $\C^n$ containing $z$, there exists an open neighborhood $U$ of $z$ in $U'$ such that  for every increasing sequence of psh functions $(u_j)_j$ on $U$ 
with $u_j \le 0$ on $K \cap U$, then 
$$(\sup_j u_j)^* \le 0$$
on $K \cap U$.
\end{definition}

Observe that $K$ is locally regular if for every $z \in K$, for every ball $B(z,r)$ centered at $z$ in $\C^n$, then $V_{K\cap B(z,r)}$ is continuous. Moreover, we have the following property.

\begin{lemma} \label{le-dktuongduonglocallyregular} The following statements are equivalent:

(i) $K$ is a locally regular set,

(ii) $(K,Q)$ is regular for every continuous function $Q$ on $K$,

(iii) $(K,-\phi_{FS})$ is regular for $\phi_{FS}(z):= -\frac{1}{2} \log(1+|z|^2)$.
\end{lemma}

The implication (iii) to (i) was proved in \cite[Theorem A]{Sadullaev-pluriregular-global}. The equivalence between (i) and (ii) was proved in  \cite[Proposition 6.1]{Dieu-regularityset}. We also refer to Theorem \ref{the-dktuongduonglocallyregular2} below for a more general result in the global setting.  

We note however that there is an example of a compact set $K$ in $\C^n$ such that $K$ is not locally regular but $(K,Q)$ for $Q \equiv 0$ is regular; see \cite[Proposition 8.1]{Sadullaev-regularset}.  One can consult \cite[Section 5]{NCNguyen-envelope} for a survey of examples of locally regular sets.

We now consider a more general setting. 
Let $(X,\omega)$ be a compact K\"ahler manifold of dimension $n$. A subset $E$ in $X$ is said to be pluripolar if $E \subset \{u=-\infty\}$ for some quasi-psh function $u$ on $X$. It is well-known that the pluripolarity is equivalent to the local pluripolarity, see \cite{Vu_pluripolar} and references therein. Let $K\subset X$ be a non-pluripolar compact set. We denote by $PSH(X, \omega)$ the set of $\omega$-psh functions on $X$. Let $\phi$ be a real continuous function on $K$.   % (or even of Hausdorff's dimension at most $2n-2$, CITATION!!!). 
Define
\begin{equation}\label{fik}
\phi_K:= \sup\{\psi \in PSH(X, \omega): \, \psi \le \phi \: \text{ on } K\}.
\end{equation}
Since $K$ is non-pluripolar, the function $\phi_K^*$ is a bounded 
$\omega$-psh function. As in Lemma \ref{le-lowersemi}, the function $\phi_K$ is lower semi-continuous.    %If $\phi_K= \phi^*_K$, then we say $(K,\phi)$ is \emph{regular}. A stronger notion is the following: we say that $K$ is \emph{locally regular} if for every $z \in K$ there is an open neighborhood $U$ of $z$ such that for every increasing sequence of psh functions $(u_j)_j$ on $U$ with $u_j \le 0$ on $K \cap U$, then $$(\sup_j u_j)^* \le 0$$on $K \cap U$. We now consider for the moment a very closely related construction in the special case when $X= \P^n$. 

\begin{remark} \label{re-sosanhphiKVK} Recall that $\phi_{FS}(z)= -\frac{1}{2} \log(1+|z|^2)$, $z \in \C^n$. Let $X= \P^n$ and $\omega:= \omega_{FS}$ the Fubini-Study form on $\P^n$.  Recall that $\ddc \phi_{FS}= -\omega_{FS}$ and  if $u \in \mathcal{L}(\C^n)$, then $u+ \phi_{FS}$ belongs to $PSH(\P^n, \omega_{FS})$. It follows that 
$$(\phi_{FS})_K= V_K+ \phi_{FS}$$  
on $\C^n$. More generally, for every $\phi$, one has
\begin{align}\label{eq-sosanhjphiKVK}
\phi_K= V_{K, \phi- \phi_{FS}}+ \phi_{FS}
\end{align}
on $\C^n$. %It was proved in \cite{Sadullaev-pluriregular-global} that for every non-pluripolar compact set $K$ in $\C^n$, there holds: $K$ is locally regular if and only if $(\phi_1)_K$ is continuous for $\phi_1 \equiv 0$ on $K$, which in turn is equivalent to the fact that $V_{K, -\phi_{FS}}$ is continuous (by (\ref{eq-sosanhjphiKVK})). We refer to \cite[Theorem 1.2]{NCNguyen-envelope} for generalizations. 
\end{remark}

\begin{definition} \label{def-localregular2} We say that $(K,Q,\omega)$ is \emph{regular} if $\phi_{K}= \phi^*_{K}$ (or equivalently, $\phi_{K}$ is continuous).  The set $K$ is said to be \emph{locally regular} if $K\cap U'$ is locally regular for every local chart $U'$ of $X$.  
%if for every $z\in K$,  every open subset $U'$ in $X$ contanining $z$, there exists an open neighborhood $U$ of $z$ in $U$ such that  for every increasing sequence of psh functions $(u_j)_j$ on $U$ 
%with $u_j \le 0$ on $K \cap U$, then 
%$$(\sup_j u_j)^* \le 0$$
%on $K \cap U$.
\end{definition}

We will drop the notation $\omega$ from $(K,Q,\omega)$ when $\omega$ is clear from the context.  As in the local case,  we have the following property proved recently by Nguyen \cite{NCNguyen-envelope}.

\begin{theorem}[{\cite[Theorem 1.2]{NCNguyen-envelope}}]  \label{the-dktuongduonglocallyregular2}
The following statements are equivalent:

(i) $K$ is a locally regular set,

(ii) $(K,Q,\omega)$ is regular for every continuous function $Q$ on $K$,

(iii) $(K,0,\omega)$ (where the weight $Q \equiv 0$) is regular.
\end{theorem}

We note that if $X= \P^n$, $\omega= \omega_{FS}$,  and $K$ is a non-pluripolar compact subset in $\C^n$, then the condition (iii) in Theorem \ref{the-dktuongduonglocallyregular2} is equivalent to the fact  that $(K,-\phi_{FS})$ is regular by Remark \ref{re-sosanhphiKVK}.

%\subsection{Local regularity}

We describe now an important class of non-pluripolar sets which are relevant in applications.
In what follows, by a (convex) \emph{polyhedron} in $\R^M$, 
we mean a subset in $\R^M$ which is the intersection of a 
finite number of closed half-hyperplanes in $\R^M$.

\begin{definition}\label{def_piecewisesmooth}
A subset $K$ of a real $M$-dimensional smooth manifold $Y$ is called a nondegenerate piecewise-smooth submanifold of dimension $m$ of $Y$ if for every point $p \in K,$ there exists a local chart $(W_p,\Psi)$ of $Y$ such that $\Psi$ is a smooth diffeomorphism from $W_p$ to the unit ball of $\R^{M}$ and  $\Psi(K \cap W_p)$ is the intersection with the unit ball  of the union of finitely many  convex polyhedra of the same dimension $m.$ %In particular, when $\Psi(K \cap W_p)$ is the intersection with the unit ball of a convex polyhedron of dimension $m$ or of the complementary polyhedron of a convex one of dimension $m$, we say that $K$ is a  $\mathcal{C}^5$ submanifold of dimension $m$ with nondegenerate piecewise boundary.
\end{definition}

A point $p \in K$ is said to be \emph{a regular point} 
of $K$ if the above local chart $(W_p, \Psi_p)$ can be 
chosen such that $\Psi_p(K\cap W_p)$ is the intersection 
of the unit ball with an $m$-dimensional vector subspace in $\R^M$, 
in other words, $K$ is a $m$-dimensional submanifold locally near $p$. 
\emph{The regular part} of $K$ is, by definition, the set of regular points of $K.$  
\emph{The singular part} of $K$ is the complement of the regular part 
of $K$ in $K.$ Hence if $K$ is a smooth manifold with boundary, 
then the boundary of $K$ is the singular part of $K$ and its complement 
in $K$ is the regular part of $K$. 

Now let $Y$ be a complex manifold of dimension $n$
and let $K$ be a nondegenerate piecewise-smooth
submanifold of $Y.$ Since $Y$ is a complex manifold,
its real tangent spaces have a natural complex structure $J$. 
We say that  $K$ is \emph{(Cauchy-Riemann) generic} 
in the sense of Cauchy-Riemann geometry if for every $p \in K$ 
and every sequence of regular points $(p_m)_m \subset K$
approaching to $p$, any limit space of the sequence of tangent spaces
of $K$ at $p_m$ is not contained  in a complex hyperplane
of the (real) tangent space at $p$ of $Y$ (equivalently, 
if $E$ is a limit space of the sequence $(T_{p_m}K)_{m\in \N}$
of tangent spaces at $p_m$, then  we have $E+ J E= T_p Y$, 
where $T_p Y$ is the real tangent space of $Y$ at $p$). 

For a generic $K$, note that the space $T_pK \cap J T_p K$ 
($p$ is a regular point in $K$) is invariant under $J$ and hence 
has a complex structure induced by $J$. In this case, 
the complex dimension of $T_pK \cap J T_p K$ is the same
for every $p$ and is called \emph{the CR dimension} of $K$.
If $r$ denotes the CR dimension of $K$, then $r= \dim K -n$.
Thus the dimension of a generic $K$ is at least $n.$ 

If $K$ is generic and  $\dim K=n$, then $K$ is said to be (maximally)
\emph{totally real}, and it is locally the graph of a smooth function over 
a small ball centered at $0 \in \R^n$ which is tangent at $0$ to $\R^n$. 
Examples of piecewise-smooth CR generic submanifolds are  polygons in 
$\C$ or boundaries of polygons in $\C$, and polyhedra of dimension 
$n$ in $\R^n \subset \C^n$.

The following result answers the question 
raised in \cite[Remark 1.8]{BoucksomBermanWitt}.

\begin{theorem}[{\cite[Theorem 2.3]{Marinescu-Vu}}] \label{the-BMpropertylocallyregu} 
Every compact nondegenerate piecewise-smooth 
(Cauchy-Riemann) generic submanifold of $X$ is locally regular. 
\end{theorem}

%Note that Theorem \ref{the-BMpropertylocallyreguCn} is a direct consequence of the above result. 
It was shown in \cite[Corollary 1.7]{BoucksomBermanWitt} 
that $K$ is locally regular if $K$ is smooth real analytic. We refer to \cite[Theorem 1.1]{VANguyenlocalregular} for the case where $K$ is smooth. 
%Theorem \ref{the-BMproperty} is actually a direct consequence of Theorem \ref{the-BMpropertylocallyregu} and the criterion \cite[Proposition 3.4]{Bloom-Levenberg-PW} giving a sufficient condition for measures to be Bernstein-Markov. 
The proof of Theorem \ref{the-BMpropertylocallyregu} presented in 
\cite{Marinescu-Vu} relies heavily on the construction of a suitable family 
of analytic discs partly attached to 
$K$, as outlined in \cite{Vu_MA,Vu_feketepoint}. 
For a comprehensive account of analytic discs, we refer to \cite{MerkerPorten,MerkerPorten2}.
%We refer to \cite[Section 4.1]{NCNguyen-envelope} for more examples of locally regular sets.
We recall now a sufficient criterion for measures to be Bernstein-Markov. 

\begin{proposition} \label{pro-tieuchuanBM} (\cite[Proposition 3.6]{Bloom-Levenberg-PW}) Let $L$ be a positive line bundle on $X$.  Assume that $K$ is locally regular and there exist constants $r_0>0$, $\tau>0$ such that $\mu(\B(z,r)) \ge r^\tau$ for every $r \le r_0$ and $x\in K$. Then $\mu$ is Bernstein-Markov with respect to $(K,\phi,L)$ for every continuous weight $\phi$.  
\end{proposition}

We refer to \cite[Proposition 3.4]{Bloom-Levenberg-PW} for a somewhat stronger criterion. The following direct consequence of 
Theorem \ref{the-BMpropertylocallyregu} and Proposition \ref{pro-tieuchuanBM} provides us a large class of Berstein-Markov measures. %the criterion \cite[Proposition 3.4]{Bloom-Levenberg-PW} 

\begin{theorem}[{\cite[Theorem 2.1]{Marinescu-Vu}}]   \label{the-BMproperty} 
Let $K$ be a compact nondegenerate piecewise-smooth generic  submanifold  of $X$. 
Then for  every continuous function $\phi$ on $K$, if $\mu$ 
is a finite measure whose support is equal to $K$ such that there 
exist constants $\tau>0, r_0>0$ satisfying 
$\mu(\B(z,r) \cap K) \ge r^\tau$ 
for every $z \in K$, and every $r \le r_0$  (where $\B(z,r)$ 
denotes the ball of radius $r$ centered at $z$ induced by 
a fixed smooth Riemannian metric on $X$), 
then $\mu$ is a Bernstein-Markov measure with respect to $(K,\phi,L)$.
\end{theorem}

Consequently, we obtain the following:

\begin{theorem}[{\cite[Theorem 1.3]{Marinescu-Vu}}] \label{the-BMpropertycn} 
Let $K$ be a compact nondegenerate  piecewise-smooth generic submanifold in $\C^n$. 
Let $\mu$ be a finite measure supported on $K$ such that there exist 
constants $\tau>0, r_0>0$ satisfying $\mu(\B(z,r)\cap K) \ge r^\tau$ 
for every $z \in K$,  $r \le r_0$  (where $\B(z,r)$ denotes the ball 
of radius $r$ centered at $z$ in $\C^n$). 
Then  for  every continuous function $Q$ on $K$, 
$\mu$ is a Bernstein-Markov measure with respect to $(K,Q)$.
\end{theorem}

Let $h_0, h, \phi$ be as in the previous section. 
Recall that the Chern form of $h_0$ is equal to $\omega$.  Define 
$$\phi_{K,k}:=\sup\big\{ k^{-1} \log |\sigma|_{h_0}: 
\quad  \sigma \in H^0(X, L^k), \,  \sup_K (|\sigma|_{h_0^k} e^{-k\phi}) \le 1\big\}.$$
Clearly $\phi_{K,k} \le \phi_K$. 
%We recall the well-known fact that the extremal functions on $X$ can be defined using sections of an ample line bundle on $X$. This generalizes a similar classical fact in $\mathbb{C}^n$ for $V_{K,Q}$. 
Here is a well-known observation, see, e.g., \cite[Lemma 3.4]{Marinescu-Vu}.

\begin{lemma} \label{le-pointwiseconvergence} 
The sequence $(\phi_{K,k})_k$ converges pointwise to $\phi_K$ as 
$k \to \infty$ and there holds $\phi_{K,k} \le \phi_K$. 
In particular,  $\phi_K$ is lower semi-continuous.
\end{lemma}

\begin{lemma}[{\cite[Lemma 3.2]{Bloom-Shiffman}} 
or {\cite[Lemma 3.5]{Marinescu-Vu}}]  \label{le-hoitudeuenvelop} 
Assume that $(K,\phi)$ is regular. 
Then $\phi_{K,k}$ converges uniformly to $\phi_K$ as $k \to \infty$. 
\end{lemma}

%Let $\bold{1}_K$ be the characteristic function of $K$. 
We put 
$$\tilde{\phi}_{K,k}:= \frac{1}{2k} \log \widetilde{B}_k=
\frac{1}{2k} \log \sum_{j=1}^{d_k} |s_j|^2_{h_0^k}.$$

%The following well-known result is our starting point.
%We recall the following well-known limit. 

\begin{proposition}[{\cite[Proposition 3.6]{Marinescu-Vu}}]\label{pro-uniforBergmann}  
Assume that  $(K,\mu,\phi)$ 
satisfies the Bernstein-Markov property. Then we have
$$\tilde{\phi}_{K,k} \ge \phi_{K,k} +O(1/k)$$
and  for every $\varepsilon>0$, 
\begin{align}\label{inetinhtangcuaphinga}
\tilde{\phi}_{K,k} \le \phi_K + \varepsilon,
\end{align}
for $k$ sufficiently large. In particular, $\tilde{\phi}_{K,k}$ converges pointwise to $\phi_K$ as $k\to \infty$ and 
\begin{equation}\label{eq:B_k}
\lim_{k \to \infty} \widetilde{B}_k^{1/k}= e^{2\phi_K}.
\end{equation}
Assume furthermore that $(K,\phi)$ is regular. Then there holds
\begin{equation}\label{eq:fiKk}
\big\| \tilde{\phi}_{K,k} - \phi_K\big\|_{\cali{C}^0(X)} \to 0,
\:\:\text{as $k \to \infty$}.
\end{equation}
\end{proposition}

%One can see from (\ref{inetinhtangcuaphinga}) that the sequence $\tilde{\phi}_{K,k}$ almost increases to $\phi_K$. 
Note that the limit in \eqref{eq:B_k} 
%last limit in Proposition \ref{pro-uniforBergmann} 
is independent of $\mu$.  
We refer to \cite[Lemma 2.8]{Beckermann-Putinar-Saff-Sty} 
for more informations in the case $K \subset \C^n \subset X =\P^n$.

%\begin{remark} \label{re-notlowerbound} Recall  $\phi_K \le \phi$ on $K$. If $x \in K$ is a point so that $\phi_K(x)< \phi(x)$, then by  Proposition \ref{pro-uniforBergmann} we see that $B_k^{-1}(x)$ grows exponentially as $k \to \infty$. Consider now the case where $K=X$ and $\phi$ is not an $\omega$-psh function. In this case there exists $x \in X$ with $\phi_X(x)< \phi(x)$, and hence $B_k$ becomes exponentially small as $k \to \infty$.   \end{remark}

\subsection{H\"older regularity of extremal psh envelopes}

Let $\alpha \in (0,1]$ and $Y$ be a metric space. For every $f: Y \to \C$, we recall that 
$$\|f\|_{\mathcal{C}^{0,\alpha}}= \|f\|_{\mathcal{C}^0}+  \sup_{x,y \in Y, x \not =y} \frac{|f(x)-f(y)|}{d(x,y)^\alpha},$$
where $d(\cdot,\cdot)$ denotes the metric on $Y$. We denote by $\mathcal{C}^{0,\alpha}(Y)$ the space of functions on $Y$ of finite $\mathcal{C}^{0,\alpha}$-norm. If $0<\alpha<1$, then for simplicity we will sometimes write $\mathcal{C}^\alpha$ for $\mathcal{C}^{0,\alpha}$. The following notion introduced in \cite{DMN} will play a crucial role for us. 

\begin{definition} Let $\alpha \in (0,1]$ and $\alpha' \in (0,1]$. A  non-pluripolar compact subset $K$ in $X$ is said to be $(\mathcal{C}^{0,\alpha}, \mathcal{C}^{0,\alpha'})$-regular if %for all H\"older continuous function $\phi$ of order $\alpha$ on $K,$  we have $\phi_K= \phi_K^*$ and 
for any positive constant $C,$ the set $\{\phi_K: \phi \in \mathcal{C}^{0,\alpha}(K) \text{ and } \|\phi\|_{\mathcal{C}^{0,\alpha}(K)} \le C\}$ is a bounded subset of $\mathcal{C}^{0,\alpha'}(X).$
%$\|\phi_K\|_{\mathcal{C}^{\alpha'}(X)} \le C \|\phi\|_{\mathcal{C}^{\alpha}(K)}+C.$ 
\end{definition}

The following result provides examples for the above regularity notion.

\begin{theorem}[{\cite[Theorem 2.3]{Vu_feketepoint}}]
\label{the1CalphaClapharegu} 
Let $\alpha \in (0,1).$ Then any  compact generic 
nondegenerate piecewise-smooth submanifold $K$ 
of $X$ is $(\mathcal{C}^{0,\alpha}, \mathcal{C}^{0,\alpha/2})$-regular.  
Moreover if $K$ is smooth, then  $K$ is 
$(\mathcal{C}^{0,\alpha}, \mathcal{C}^{0,\alpha})$-regular.
\end{theorem}

If $K=X$, and $\phi \in \mathcal{C}^{0,1}$, then it was shown in 
\cite{Chu-Zhou-optimal-envelope} that $\phi_X \in \mathcal{C}^{0,1}$, 
hence $X$ is $(\mathcal{C}^{0,1}, \mathcal{C}^{0,1})$-regular; 
see also \cite{Berman-C11regula,Chu-Zhou-optimal-envelope,Tosatti-envelop}
for more information. In the case where $K=X$ or $K$ is an open subset 
with smooth boundary in $X$ it was proved in \cite{DMN} 
that $K$ is $(\mathcal{C}^{0,\alpha}, \mathcal{C}^{0,\alpha})$-regular 
for $\alpha \in (0,1)$, see also \cite{Lu-To-Phung}. A necessary and sufficient condition for $(\mathcal{C}^{0,\alpha}, \mathcal{C}^{0,\alpha'})$-regularity  was obtained in \cite[Theorem 1.2]{NCNguyen-envelope}. Thanks to this criterion, one has the following new examples for regular sets. 

\begin{theorem} [{\cite[Theorem 1.4]{NCNguyen-envelope}} and 
{\cite[Theorem 3.2]{AhnNguyen}}] 
Let $K$ be a compact uniformly polynomially cuspidal set in $\mathbb{R}^n$.
Then, $K$ is $(\mathcal{C}^{0,\alpha}, \mathcal{C}^{0,\alpha'})$-regular 
(as a subset of $\mathbb{C}^n \subset \mathbb{P}^n$) 
for some constants $\alpha \in (0,1]$ and $\alpha' \in (0,1]$. 
\end{theorem}

 We recall  that bounded convex sets with non-empty interior 
 and fat subanalytic sets are uniformly cuspidal 
 (see \cite[Corollary 6.6]{PawluckiPle_markov} and 
 \cite[Corollary 3.1]{Plesniak-dissertation}).  
 We refer to \cite{Bierstone_Milman} and \cite{PawluckiPle_markov} 
 for more information on subanalytic sets and uniformly polynomialy cuspidal ones respectively.

It is not known if 
Theorem \ref{the1CalphaClapharegu} holds for $\alpha=1$. 
However, we have the following.
%===
\begin{theorem}[{\cite[Theorem 3.11]{Marinescu-Vu}}] 
\label{the1CalphaClapharegunanchinhquy} 
Let $\delta \in (0,1), C_1>0$ be constants. Then the following properties hold: 
Let $K$ be a compact generic non-degenerate piecewise-smooth submanifold
of $X$. Then there exists a constant $C_2>0$
such that for every $\phi \in \mathcal{C}^{1, \delta}(K)$ with
$\|\phi\|_{\mathcal{C}^{1,\delta}} \le C_1$, then 
$\|\phi_K\|_{\mathcal{C}^{1/2}} \le C_2$. Moreover, if $K$ is smooth, then $\phi_K \in \mathcal{C}^{0,1}$ with $\|\phi_K\|_{\mathcal{C}^{0,1}} \le C_2$.    
%In particular, if $K$ is a compact nondegenerate piecewise-smooth generic submanifold in $\C^n$.Then $V_K \in \mathcal{C}^{1/2}(\C^n)$. Additionally if $K$ is smooth, then $V_K \in \mathcal{C}^{0,1}(\C^n)$. 
\end{theorem}
%===

We note that the assertion concerning the estimate $\|\phi_K\|_{\mathcal{C}^{1/2}} \le C_2$ was not stated in \cite[Theorem 3.11]{Marinescu-Vu} but this is deduced by the same arguments given in the paragraph before \cite[Theorem 3.11]{Marinescu-Vu}. Indeed, the only reason why \cite[Theorem 2.3]{Vu_feketepoint} was stated for $\alpha \in (0,1)$   because its proof uses \cite[Lemma 2.6]{Vu_feketepoint} about the regularity of subharmonic functions on the unit disk in $\C$, which only works for $\alpha \in (0,1)$. However, if the regularity of the weight $\phi$ increases sufficiently (i.e., in $\mathcal{C}^{1,\delta}$ for some $\delta>0$), then an improvement of \cite[Lemma 2.6]{Vu_feketepoint}   holds (see \cite[Lemma 3.3]{Marinescu-Vu}). This allows us to treat the case $\alpha=1$.

At this juncture, we mention an example from \cite{Sadullaev-example} 
of a domain $K$ with a $\mathcal{C}^0$ boundary, yet $(K,\phi)$ 
fails to be regular, even for $\phi = 0$. 
Applying Theorem \ref{the1CalphaClapharegunanchinhquy}
to $X = \P^n$, we derive the following result, which directly 
implies \cite[Conjecture 6.2]{Sadullaev-Zeriahi} as a special case. 

\begin{theorem} \label{the1CalphaClaphareguCn} (\cite[Theorem 3.12]{Marinescu-Vu}) 
Let $K$ be a compact generic nondegenerate  piecewise-smooth submanifold in $\C^n$.
Then $V_K \in \mathcal{C}^{1/2}(\C^n)$. 
Additionally if $K$ is smooth, then $V_K \in \mathcal{C}^{0,1}(\C^n)$. 
\end{theorem}

Recall that $$V_K= \sup \{\psi \in \mathcal{L}(\C^n): \psi \le 0 \, \text{ on } K\}.$$
We note that the fact that $V_K \in \mathcal{C}^{0,1}(\C^n)$ when $K$ is smooth was proved in \cite{Sadullaev-Zeriahi} and \cite{Vu_feketepoint} independently.
%We mention at this point an example from \cite{Sadullaev-example}
%of a domain $K$ with $\mathcal{C}^0$ boundary but $(K,\phi)$ is not regular even 
%for $\phi =0$.  Applying Theorem \ref{the1CalphaClapharegunanchinhquy} 
%to $X= \P^n$ we obtain the following result that 
%implies \cite[Conjecture 6.2]{Sadullaev-Zeriahi} as a special case. 

%We give now a Bernstein-Markov type inequality. 
%We will not use it anywhere in the paper.

\subsection{Bernstein-Markov inequality}

Bernstein-Markov type inequalities give an upper bound for the derivatives of polynomials in terms of their supnorms. They are objects of great interest in 
approximation theory. There is a large literature on this topic, e.\,g.\ 
\cite{Berman-OrtegaCerda,Bos-Levenberg-Milman-Taylor2,
Bos-Levenberg-Milman-Taylor,
Brudnyi,Brudnyi-bernstein-ine-sub,
Coman-Poletsky,Pierzchala-holder,Yomdin-smooth-parametrization}, 
to cite just a few. We also refer to \cite{Totik-survey} for a survey about Bernstein-Markov inequalities in dimension one. We are interested in higher dimensional versions of Bernstein-Markov inequalities.

\begin{theorem}[Bernstein-Markov type inequality, {\cite[Theorem 3.13]{Marinescu-Vu}}]
\label{the-Bernstein-totally-real}
Let $K$ be a compact generic nondegenerate piecewise-smooth submanifold in $\C^n$. 
Then there exists a constant $C>0$ such that for every 
complex polynomial $p$ on $\C^n$ we have
\begin{equation}\label{eq:grad1}
\|\nabla p\|_{L^\infty(K)} \le C (\deg p)^{2}\|p\|_{L^\infty(K)}.
\end{equation} 
If additionally $K$ is smooth, then 
\begin{equation}\label{eq:grad2}
\|\nabla p\|_{L^\infty(K)} \le C \deg p \,\|p\|_{L^\infty(K)}.
\end{equation}  
\end{theorem} 

Note that the exponent of $\deg p$ is optimal as it is well-known for
the classical Markov and Bernstein inequalities in dimension one. 
The above result was known when $K$ is algebraic in $\R^n$, 
see \cite{Berman-OrtegaCerda,Bos-Levenberg-Milman-Taylor}. 
We emphasize that similar estimates also hold for other situations 
(with the same proof), for example, $K= \S^{n-1} \subset \R^n$ by considering $K$ as a maximally totally 
real submanifold in the complexification of $\S^{n-1}$.

\proof  %Let $$\phi_K:= \sup\{\psi\in\mathcal{L}(\C^n):\sup_K \psi \le 0\}$$which is
Recall that $$V_K= \sup \{\psi \in \mathcal{L}(\C^n): \psi \le 0 \, \text{ on } K\}.$$
Observe first that $V_K$ is $\mathcal{C}^{0,1}$ if $K$ is smooth or
$\mathcal{C}^{1/2}$ in general  by Theorem \ref{the1CalphaClaphareguCn}.
Let $p$ be a complex polynomial in $\C^n$. Put $k:= \deg p$. 
Since $\frac{1}{k} (\log |p|- \log \max_K|p|)$ is a candidate 
in the envelope defining $V_K$, we get
$$|p| \le e^{k V_K} \max_K |p|$$
on $\C^n$. 
We use the same notation $C$ to denote a constant depending only on $K,n$.  
Let $a=(a_1,\ldots,a_n) \in K \subset \C^n$. Let $r>0$ be a small constant.  
Consider the analytic disc $D_a:= (a_1+ r\D, a_2,\ldots, a_n)$. 
Applying the Cauchy formula to the restriction of $p$ to $D_a$ shows that 
$$|\partial_{z_1}p(a)| \le r^{-1} \max_{D_a} |p|\le
r^{-1} (\max_K |p|) \max_{D_a} e^{k V_K}.$$
Since $\phi_K = 0$ on $K$, using $\mathcal{C}^{1/2}$ regularity
of $V_K$ gives
$$|\partial_{z_1}p(a)| \le r^{-1} \max_{D_a} |p|\le 
r^{-1} (\max_K |p|) e^{C k r^{1/2}}$$
for some constant $C>0$ independent of $p$ and $a$. 
Choosing $r= k^{-2}$ in the last inequality yields
$$|\partial_{z_1}p(a)| \le C k^{2} \max_{D_a} |p|.$$
Similarly we also get  $$|\partial_{z_j}p(a)| \le 
C k^{2} \max_{D_a} |p|$$
for every $1 \le j \le n$. Hence the first desired inequality
\eqref{eq:grad1} 
for general $K$. When $K$ is smooth, the arguments are similar. 
This finishes the proof.  
\endproof

It should be noted that Theorem \ref{the-Bernstein-totally-real} 
only concerns complex polynomials. If we consider real polynomials in 
$\R^{2n} \approx \C^n$, then an inequality of type \eqref{eq:grad2})
can be used to characterize the real algebraicity of $K$ as showed by
Bos-Levenberg-Milman-Taylor \cite{Bos-Levenberg-Milman-Taylor} in the following result. 

\begin{theorem} [{\cite{Bos-Levenberg-Milman-Taylor}}] 
 Let $K$ be a compact smooth submanifold of $\R^N$. 
 Then the following statements are equivalent:

(i) $K$ is algebraic.

(ii) $K$ satisifies the following Bernstein-Markov inequality: 
There exists a constant $C>0$ such that for every (real) polynomial $p$ on $\R^N$, we have 
  \begin{equation}\label{eq:grad3}
\|\nabla p\|_{L^\infty(K)} \le C \deg p \,\|p\|_{L^\infty(K)}.
\end{equation}  
\end{theorem}

A stronger version for not necessarily compact manifolds was proved in \cite{Bos-Levenberg-Milman-Taylor2}.

 \subsection{Polynomial growth of Bergman kernel functions}
\label{subsec-polygrowth}

%We commence by exhibiting a large class of Bernstein-Markov measures. 
As established by Bloom and Levenberg \cite{BlLe15}, for Bernstein-Markov measures, 
the equidistribution of zeros of random sections in probability holds.
To derive a rate of convergence,
we will need finer properties of Bergman (Christoffel-Darboux) kernel.

%\begin{theorem} \label{the-BMproperty}   
%Let $K$ be a compact nondegenerate smooth piecewise-smooth 
%Cauchy-Riemann generic   submanifold  of $X$. 
%Then for  every continuous function $\phi$ on $K$, if $\mu$ 
%is a finite measure whose support is equal to $K$ such that there 
%exist constants $\tau>0, r_0>0$ satisfying 
%$\mu(\B(z,r) \cap K) \ge r^\tau$ 
%for every $z \in K$, and every $r \le r_0$  (where $\B(z,r)$ 
%denotes the ball of radius $r$ centered at $z$ induced by 
%a fixed smooth Riemannian metric on $X$), 
%then $\mu$ is a Bernstein-Markov measure with respect to $(K,\phi,L)$.
%\end{theorem}
%===
The Monge-Amp\`ere current $(\ddc \phi_K^* + \omega)^n$ is called 
\emph{the equilibrium measure} associated to $(K,\phi)$. 
It is well-known that this measure is supported on $K$.
By \cite[Theorem B]{BoucksomBermanWitt} one has 
\begin{align}\label{conver-weak-bergmankernelintro}
d_k^{-1} B_k \mu \to (\ddc \phi_K^*+ \omega)^n\,,\quad k\to\infty,
\end{align}
provided that $\mu$ is a Bernstein-Markov measure associated to $(K,\phi,L)$.
The last property suggests that the Bergman kernel function $B_k$ 
cannot behave too wildly as $k\to\infty$.  
%It is then a natural question to understand asymptotic properties 
%of $B_k$ or at least obtain upper 
%and lower bounds for $B_k$. In this direction,  
%Our next main result shows that a similar property 
%holds if $X$ is replaced by a generic CR nondegenerate 
%smooth  pieceweise smooth submanifold. 

\begin{theorem} \label{the-strongBMintro}  
Let $K$ be a compact nondegenerate piecewise-smooth 
Cauchy-Riemann generic   submanifold of $X$. Let $n_K$ be the dimension of $K$. 
Let $\phi$ be a H\"older continuous function of H\"older exponent $\alpha \in (0,1)$ on $K$, 
let $\Leb_K$ be a smooth volume form on $K$, and  $\mu = \rho \Leb_K$, 
where $\rho \ge 0$ and $\rho^{-\lambda} \in L^1(\Leb_K)$ for some
constant $\lambda>0$. Then, there exists a constant $C>0$ such that 
$$\sup_K B_k \le C k^{2n_K (\lambda+1)/(\alpha\lambda)}$$
for every $k$. 
\end{theorem}

%When $\phi$ is Lipschitz (in particular when $\phi \equiv 0$), we have much more precise bounds. 

Note that by  the  proof of \cite[Theorem 1.3]{Dinh-Ma-Marinescu} 
or \cite[Theorem 3.6]{DMN}, for every H\"older continuous function 
$\phi_1$ on $X$, and  $\mu_1:= \omega^n$,  the Bergman kernel function 
of order $k$ associated to  $(X, \mu_1, \phi_1)$ grows at most polynomially 
on $K$ as $k \to \infty$; see also \cite[Theorem 3.1]{BoucksomBerman} 
for the case where $\phi_1$ is smooth.

\begin{theorem} \label{the-bergman-smooth} 
Assume that the following two conditions hold:

(i) $K$ is maximally totally real and smooth,

(ii) $\phi \in \mathcal{C}^{1,\delta}(K)$ for some constant $\delta>0$.

\noindent
Then there exists $C>0$ such that  for every $k$ and every $x \in K$
the following holds 
$$B_k(x) \le C k^{n}\,.$$
\end{theorem}

Consider the case when $X= \P^n$, $L:=\mathcal{O}(1)$, $h_0=h_{FS}$ 
is the Fubini-Study metric on $\mathcal{O}(1)$, 
and $K$ is a smooth maximally totally real compact submanifold in 
$ \C^n \subset \P^n$, $Q$ is a continuous function on $K$ and $h:= e^{-2 \phi} h_0$ on $K$, where  $\phi:= Q -\frac{1}{2}\log (1+|z|^2)$. Observe that $\phi$ is
in $\mathcal{C}^{1,\delta}(K)$ if $Q$ is so. 
In this case the hypothesis of Theorem  \ref{the-bergman-smooth} are fulfilled.  
%Thus Theorem \ref{the-bergman-smooth}implies Theorem \ref{the-bergman-smoothCn}.

Here is a quantitative version of Lemma \ref{le-pointwiseconvergence}.
 
\begin{proposition} \label{pro-hoituC0phiKkphiK}  
Let $K$ be a compact generic nondegenerate
piecewise-smooth submanifold of $X$.  Let $\phi$ be a 
H\"older continuous function on $K$. Then, we have 
$$\big\| \phi_{K,k}- \phi_K\big\|_{\cali{C}^0(X)}=
O\left(\frac{\log k}{k}\right).$$
\end{proposition}
 The above result was proved for $K=X$ in 
\cite[Corollary 4.4]{Dinh-Ma-Marinescu}.

\begin{theorem}[{\cite[Theorem 2.6]{Marinescu-Vu}}] 
\label{the-convergencebergmankernel}  
Let $K$ be a compact  nondegenerate piecewise-smooth 
generic submanifold of $X$.   Let $\phi$ be a H\"older 
continuous function on $K$. Let $\mu$ be a smooth volume form on $K$. 
Then we have 
$$\left\| \frac{1}{2k} \log \widetilde{B}_k - \phi_K \right\|_{\mathcal{C}^0(X)}=
O\left(\frac{\log k}{k}\right),$$
as $k \to \infty$, where 
$$\widetilde{B}_k:= e^{2 k \phi} B_k=
\sup\big\{ |s(x)|_{h_0^k}^2: \: s \in H^0(X,L^k), \: \|s\|_{L^2(\mu,h^k)}=1\big\}.$$
\end{theorem}  
 
The proof of this result is based essentially on the polynomial upper bound for $B_k$ in Theorem \ref{the-strongBMintro}.
%\begin{proof} The desired estimate is deduced directly by using Proposition \ref{pro-hoituC0phiKkphiK}, Theorem \ref{the-strongBMintro}, and following the same arguments as in the proof of Proposition \ref{pro-uniforBergmann}. We just briefly recall here how to do it. Firstly as in the proof of Proposition \ref{pro-uniforBergmann}, we have 
%$$\tilde{\phi}_{K,k}- \phi_{K,k} \ge O(k^{-1}).$$
%It remains to bound from above $\tilde{\phi}_{K,k}- \phi_{K,k}$. Combining the polynomial upper bound for $B_k$ in Theorem \ref{the-strongBMintro} and (\ref{ine-BkBergman}), one gets,  for some constants $C,N>0$ independent of $k$, 
%\begin{align*}
%\sup_K (|s|_{h_0^k} e^{- k \phi}) \le C k^{N} \|s\|_{L^2(\mu, k \phi)}
%\end{align*}
%for every $s \in H^0(X, L^k)$. This coupled with (\ref{ine-sh0wahsline}) yields
%\begin{align*}
%|s|_{h_0^k} \le k^{N} e^{k \phi_{K,k}}
%\end{align*}
%on $X$. It follows that 
%$$\tilde{\phi}_{K,k}=\frac{1}{2k}\log \widetilde{B}_k \le
%\phi_{K,k}+ N\,\frac{\log k}{k}.$$
%This finishes the proof. \end{proof}

\section{Zeros of random polynomials of several variables}
\label{S:Zrp}

Let $K$ be a non-pluripolar set in $\C^n$ and let $\mu$ 
be a probability measure on $\C^n$ such that the support of $\mu$
is contained in $K$ and is non-pluripolar. 
Let $Q$ be a continuous weight on $K$. 
\emph{We assume throughout this section that $(K,Q,\mu)$ is Bernstein-Markov.}   
Let $\mathcal{P}_k(K)$ 
be the space of restrictions of complex polynomials of degree at most $k$ 
in $\C^n$ to $K$. Let $d_k:= \dim \mathcal{P}_k(K)$, 
and let  $s_1,\ldots, s_{d_k}$ be an orthonormal basis of 
$\mathcal{P}_k(K)$ with respect to the $L^2(\mu,k Q)$-scalar product
\eqref{eq:l2muk}.  
Consider the random polynomials 
\begin{align}\label{eq-randompoly1}
p_k:= \sum_{j=1}^{d_k} a_{kj} s_j,
\end{align}
where $a_{kj}$ are complex i.i.d.\ random variables.

\subsection{Almost sure convergence}
We assume that for any $j=1,\ldots,d_k,$
the distribution of $a_{kj}$ is $f \lambda_\C$,
where $f$ is a nonnegative bounded Borel function on $\C$.   
satisfying the following mild regularity property: There exists a constant
$C>0$ such that for every $r>0$ we have,
\begin{align}\label{ine-BLdk}
\int_{\{|z|> r\}} |f| \lambda_\C \le C/r^{2}.
\end{align}   
This condition was introduced in \cite{Bloom,Bloom-Levenberg-random}. 
%We want to study the distribution of zeros of $p_k\in\mathcal{P}_k$ as $k \to \infty$. 
We denote by $[p_k=0]$ the current of integration along the zero divisor 
$\Div(p_k)$ of $p_k$. 
If $n=1$, then $[p_k=0]$ is the sum of Dirac masses at zeros of $p_k$,
counted with multiplicities.  
%For holomorphic function $f$ in $\C^n$ we let $[f=0]$ denote 
%the current of integration along the divisor  $\{f=0\}$. 

We consider the probability space $(\C^{d_k}, \mu_k)$, 
where $\mu_k$ is the probability measure with density 
$f(z_1) \cdots f(z_{d_k})$ on $\C^{d_k}$. 
Let $(\Omega, \cali{P}_\infty):= \prod_{k=1}^\infty (\C^{d_k}, \mu_k)$ 
the product of probability spaces $(\C^{d_k}, \mu_k)$. 
We want to study the behaviour of the sequence of zeros $k^{-1}[p_k=0]$ 
for a random choice of $\bold{a}=\big((a_{kj})_{j=1}^{d_k}\big)_{k \in \N}$ 
in $(\Omega, \cali{P}_\infty)$. 
The following result hightlights the importance of the notion 
of Bernstein-Markov measures. 

\begin{theorem}[{\cite[Theorem 4.2]{Bloom-Levenberg-random}}] 
\label{the-BMBloomLevenberg}
Let $(K,Q,\mu)$ be Bernstein-Markov. 
We have  
\begin{align}\label{converg-zero}
k^{-1}[p_k=0] \to \ddc V_{K,Q}^*\,,\quad k \to \infty\,,
\end{align}
almost surely, where the convergence is the weak convergence of currents. 
In other words, for every smooth form $\Phi$ of degree $(2n-2)$
with compact support in $\C^n$, one has 
$$k^{-1}\int_{\Div(p_k)} \Phi \to 
\int_{\C^n} \ddc V_{K,Q}^* \wedge \Phi\,,\quad k \to \infty\,.$$
\end{theorem}

We refer to \cite{Bloom} for an earlier result. In dimension one, it was known that the almost sure convergence of random zeros holds if and only if 
$$\int_\C\log(1+|z|)f \Leb_\C < \infty,$$
see \cite{Ibragimov-Zaporozhetz} for the case of Kac polynomials and \cite{Dauvergne} for the general case. An optimal condition for $f$ in higher dimensions is not known. On the other hand,  it was recently proved by Bloom-Dauverge-Levenberg \cite{Bloom-Dauverge-Levenberg} that 
\begin{align}
k^{-1}[p_k=0] \to \ddc V_{K,Q}^*\,,\quad k \to \infty\,,
\end{align}
in probability provided that 
$$\int_{\{\log(1+|z|)>r\}} d \Leb_\C =o(r^{-n}).$$
This condition is optimal if $n=1$; see \cite{Bloom-Dauvergne}.

There are various (non-quantitative) generalizations of Theorem \ref{the-BMBloomLevenberg} in higher dimensions. We refer to \cite{Bayraktar-Bloom-Levenberg-2024,Bayraktar-et-al-convexbodies} for some recent developments, and to  \cite{Bayraktar-zeros-indiana,Bayraktar-mass-equi,Bayraktar-Coman-Mariescu} 
for more general setting with zeros of random holomorhic sections of a positive line bundle. 
A fundamental element in the proofs of all the aforementioned results is the asymptotic behavior of Bergman (Christoffel-Darboux) kernel functions. 
In the higher dimensional setting, quantitative estimates for the equidistribution
of zeros similar to the situation in dimension one, as presented in \cite{Shiffman-Zelditch},
are lacking. 
This can be partially seen by the fact that the asympotics of Bergman kernel function 
(for a general compact $K$) has been still far from being well understood as 
in the dimension one case.   In what follow we present 
generalizations of some results from  \cite{Shiffman-Zelditch} obtained in \cite{Marinescu-Vu}. 
Indeed, the asymptotics of the Bergman kernel function (for a general compact $K$)
as in the dimension one case remain to be fully elucidated.
We refer to Subsection \ref{subsec-polygrowth} for further details on this matter.
%In what follows, we present generalizations of certain results from \cite{Shiffman-Zelditch} obtained in \cite{Marinescu-Vu}. 

%We now discuss a large deviation type result for distribution of random zeros. 
%To state our result we need some hypothesis on $\mu$ and the distribution of the random variables $\alpha_j$.   

Our goal now is to obtain a rate of convergence in \eqref{converg-zero}.  We now recall the following notion of distance on the space of currents.
For every $\beta  \ge 0$, and $T,S$ closed positive currents 
of bi-degree $(m,m)$ on  the complex projective space $\P^n$, define
$$\dist_{-\beta}(T,S):=
\sup_{\Phi:\,\|\Phi\|_{\mathcal{C}^{[\beta], \beta- [\beta]}} \le 1}
|\langle T-S, \Phi \rangle|,$$  
where $[\beta]$ denotes the greatest integer less than 
or equal to $\beta$, and $\Phi$ is a smooth form of degree $(2n-m)$ on $\P^n$. 
It is a standard fact that the distance $\dist_{-\beta}$ for 
$\beta>0$ induces the weak topology on the space of closed positive currents 
(see for example \cite[Proposition 2.1.4]{DinhSibony_Pk_superpotential}). 
We have the following interpolation inequality: 
for every $0< \beta_1 \le \beta_2$, there exists a constant
$c_{\beta_1,\beta_2}>0$ such that
\begin{align} \label{ine-interpolar}
\dist_{-\beta_2}  \le \dist_{-\beta_1} \le c_{\beta_1,\beta_2} 
[\dist_{-\beta_2}]^{\beta_1/\beta_2};
\end{align} 
see \cite[Lemma 2.1.2]{DinhSibony_Pk_superpotential} or 
\cite{Lunardi-book-interpolation,Triebel}.

Note that the currents $[p_k=0]$ and $\ddc V_{K,Q}^*$ 
extend trivially through the hyperplane at infinity 
$\P^n \backslash \C^n$ to be closed positive currents 
of bi-degree $(1,1)$ on $\P^n$. Hence one can consider the distance $\dist_{-\beta}$ between $k^{-1}[p_k=0]$ 
and $\ddc V_{K,Q}$ as closed positive currents on $\P^n$.

\subsection{Expectation of random zeros}
Denote by $\E_k(k^{-1}[p_k=0])$ the expectation of the 
random normalized currents $k^{-1}[p_k=0]$. 
For a sequence $(S_k)_{k\ge 1}$ of currents in $\C^n$ and a sequence of positive numbers $(A_k)_k$, we write 
$$S_k=O(A_k),\quad k\to\infty,$$
 if each $S_k$ is of order $0$ and of the degree $m$ for some integer $m$,
and for every compact subset $W$ in $\C^n$,  there exists a constant $C_W>0$
such that for every smooth form $\Phi$ 
of degree $(2n-m)$ with support in $W$ and $\|\Phi\|_{\mathcal{C}^2} \le 1$, we have 
$$| \langle S_k, \Phi \rangle | \le C_W\, A_k,$$
for any $k \ge 1$.
%for some constant $C$ independent of $k, \Phi$. 
Let $\mu_k$ be the joint-distribution of $a_{k1}, \ldots, a_{kd_k}$. Hence we have $$\mu_k= f(z_1) \ldots f(z_{d_k}) \Leb_{\C^{d_k}}.$$

\begin{theorem}\label{the-expectedzeros}  
Let $(K,Q,\mu)$ be Bernstein-Markov. Assume that 

(H1) For every $k \in \N$ there exists $C_k>0$ such that 
for every $u \in \C^k$ with $\|u\|=1$,
$$\int_{\C^{d_k}}\big|\log|\langle a, u \rangle|\big|d\mu_k\le C_k.$$   

(H2) let $K$ be a non-degenerate piecewise-smooth 
generic submanifold of $\C^n$, and $Q$ be a 
H\"older continuous weight on $K$. Let $\mu= \rho \lambda_K$, 
where $\rho^{-\lambda} \in L^1(\lambda_K)$ for some constant $\lambda>0$. 

Then we have 
\begin{align} \label{eq-corexpeczeros}
\E_k\big(k^{-1}[p_k=0]\big) = 
\ddc V_{K,Q}+ O\left(\frac{C_k+\log k}{k}\right).
\end{align}
%In particular if $C_k$ can be chosen to be equal to $O(k^s)$ for some constant $s\ in (0,1)$ independent of $k$, then the error current is equal to $O(k^{^-1+s})$. 
\end{theorem}

Condition (H1) was introduced in \cite{Bayraktar-Coman-Mariescu}.   We refer to this paper for a version of Theorem \ref{the-expectedzeros} when $\C^n$ is replaced by a compact K\"ahler manifold $X$ and $K=X$. One can also find in \cite{Bayraktar-Coman-Mariescu} examples of measures satisfying (H1). In particular, by \cite[Lemma 4.8]{Bayraktar-Coman-Mariescu}, if $a_{k1},\ldots,a_{k d_k}$ have  Gaussian joint density (i.e, $\mu_k= \pi^{-d_k} e^{-\|a\|^2} \Leb_{\C^{d_k}}(a)$), then (H1) is fulfilled for $C_k=C$ independent of $k$.

Condition (H2) is a natural extension of the classical setting of 
Kac polynomials where $K$ is the unit circle in $\C$. 
In fact, in \cite{Shiffman-Zelditch} the authors considered the setting 
where $\mu$ is the surface area on a closed analytic curve in $\C$ 
that bounds a simply connected domain $\Omega$ in $\C$, or 
$\mu$ is the restriction of the Lebesgue measure on $\C$ to $\Omega$. 
This setting is relevant to random matrix theory 
as already pointed out in \cite{Shiffman-Zelditch}. 
We refer to \cite{Bloom-Dauvergne,Pritsker-Ramachandran,Pritsker-Ramachandran2} 
for partial generalizations (without quantitative estimates) 
to domains with smooth boundary in $\C$. Thus we see that Theorem \ref{the-expectedzeros}  extends \cite[Theorems 1 and 2]{Shiffman-Zelditch} to higher dimensions (except that we only obtain the error term 
$O(\frac{\log k}{k})$ instead of $O(k^{-1})$);  see also Theorem \ref{the-zerosgiaoL} below.

\proof The proof is standard once we have proved a polynomial growth for $B_k$. We recall it here for the reader's convenience. Let
$$\|p_k\|^2:= \sum_{j=1}^{d_k} |s_j|^2$$
and
$$\varphi:= k^{-1}\int_{\C^{d_k}} \log |p_k| d \mu_k =k^{-1} \int_{\C^{d_k}} \log (|p_k|/\|p_k\|) d \mu_k+ k^{-1} \log \|p_k\|.$$
By (H2) and Theorem \ref{pro-hoituC0phiKkphiK}, we see that 
 $$\|k^{-1} \log \|p_k\|- V_{K,Q}\|_{L^\infty(\C^n)} = O(\log k/k).$$
The assumption (H1) implies that
$$k^{-1} \int_{\C^{d_k}} \log (|p_k|/\|p_k\|) d \mu_k \le C_k/k.$$
It follows that 
$$\E_k\big(k^{-1}[p_k=0]\big) = \ddc \varphi =  \ddc V_{K,Q}+ O\left(\frac{C_k+\log k}{k}\right)$$
as desired.  
\endproof

Note that in the case where $a_{kj}$ are Gaussian variables
the decay rate obtained in \cite{Shiffman-Zelditch} 
is $O(k^{-1})$, and that this error term is optimal in dimension one
(this can be seen by carefully examining the calculations in
\cite[Proposition 3.3]{Shiffman-Zelditch}). 
%\begin{lemma} \label{le-distribuzmua22} Let $\delta>0$ be a constant. Assume that there is a constant $C>0$ such that  
%\begin{align} \label{ine-dkfexpo}
%|f(z)| \le C e^{-\delta |z|}
%\end{align}
% for $z \in \C$. Then for every constant $\rho \in (0,1)$, there exiss positive constant $c_1,c_2$ such that 
%$$\int_{\C^n} e^{c_1 |\log\langle a, u \rangle|^{2\rho}} d\mu_k \le a_2$$
%for every $u \in \C^{d_k}$ with $\|u\|=1$. 
%\end{lemma}

%QUESTION: can we choose $\rho=2$?
%We underline that recent results on large deviation principles for random polynomials in dimension one require the condition (\ref{ine-dkfexpo}). CITE Rapha\"el Butez ET all !!!

%\proof Let $\phi:= |\log\langle a_1 \rangle|^\rho$. Observe that $\varphi \in W^{1,2}_{loc}(\C)$. 
%\endproof

\begin{lemma} \label{le-distribuzmua2} Assume that there is a constant $C \ge 1$ such that  
\begin{align}\label{ine-dktrenfzmu3}
\int_{\{\log |z| >r\}} f d\Leb_\C \le C e^{-r}; \quad |f| \le C
\end{align}
 for  every real number $r$. Then
 (H1) is satisfied for $C_k=C_0\log k$ for a constant $C_0$ independent of $k$ and for every number $r \ge 0$ one has 
\begin{align}\label{ine-bichanmulogau}
\int_{ \{a \in \C^{d_k}:\, \log |\langle a, u \rangle| \le - r\}}  d \mu_k \le  \pi C d_k e^{-2r}
\end{align} 
for every $u \in \C^{d_k}$ with $\|u\|=1$.
\end{lemma}

\proof  The desired result is a direct consequence of \cite[Lemma 4.15]{Bayraktar-Coman-Mariescu}. Since the proof is short, we include it here for the readers' convenience.  We first check (\ref{ine-bichanmulogau}). Since $\|u\|=1$, there exists an index $1\le j \le d_k$ so that $|u_j| \ge d_k^{-1/2}$. Without loss of generality, we can assume that $j=1$. Denote by $I_r$ the left-hand side of (\ref{ine-bichanmulogau}).  Using the change of variables $a'_1:= \langle a, u \rangle$, $a'_j:= a_j$ for $j \ge 2$, we obtain 
\begin{align*}
I_r &= \int_{\C^{d_k-1}} f(a_2) \ldots f(a_{d_k-1}) d\Leb_{\C^{d_k-1}} \int_{\{a'_1 \in \C: \,\log |a'_1| \le - r\}} |u_1|^{-2} f(a'_1- \sum_{j=2}^{d_k}a_j u_j)  d \Leb_\C \\
&\le C  d_k\int_{\{a'_1:\, \log |a'_1| \le - r\}}  d \Leb_\C \\
&\le \pi C d_k e^{-2r}
\end{align*} 
by (\ref{ine-dktrenfzmu3}). Hence (\ref{ine-bichanmulogau}) follows.   Put 
$$J_r:=\int_{ \{a \in \C^{d_k}:\, \log |\langle a, u \rangle| \ge  r\}}  d \mu_k. $$
Observe that
$$\{a \in \C^{d_k}:\, \log |\langle a, u \rangle| \ge  r\} \subset \bigcup_{j=1}^{d_k} \{a: \,|a_j| \ge e^{r- \frac{1}{2}\log d_k}\}.$$
Consequently, by (\ref{ine-dktrenfzmu3}), we get
$$J_r \le d_k \int_{\{a\in \C: \,|a| \ge e^{r- \frac{1}{2}\log d_k} \}} f(a)d\Leb_\C \le C d_k e^{-r+ \frac{1}{2}\log d_k} \le C d_k^2 e^{-r}.$$
This combined with (\ref{ine-bichanmulogau}) implies
$$\int_{ \{a \in \C^{d_k}:\, \big|\log |\langle a, u \rangle|\big| \ge  r\}}  d \mu_k \le  2 \pi C d^2_k e^{-r}$$
for every number $r\ge 0$. Applying this inequality to $r=r_m:= m\log d_k $ for $m\in \N$ gives
\begin{align*}
\int_{\C^{d_k}}\big|\log |\langle a, u \rangle|\big| d\mu_k &\le \sum_{m \ge 0} \int_{ \{a \in \C^{d_k}:\,  r_m \le \big|\log |\langle a, u \rangle|\big| \le  r_{m+1}\}}\big|\log |\langle a, u \rangle|\big|  d \mu_k \\
&\le 3\log d_k+ \sum_{m \ge 3} \int_{ \{a \in \C^{d_k}:\,  r_m \le \big|\log |\langle a, u \rangle|\big| \le  r_{m+1}\}}\big|\log |\langle a, u \rangle|\big|  d \mu_k \\
&\le 3 \log d_k+  \log d_k \sum_{m \ge 3} (m+1) (2 \pi C d_k^{-m+2})\\
& \lesssim \log d_k+ d_k^{-1/2}\log d_k \lesssim \log d_k.
\end{align*}
Hence the first desired assertion follows. 
%We point out how to deduce our result from theirs. It suffices to choose $\nu=1$ and replace $R^\rho$ by $e^{-r}$ in the proof of \cite[Lemma 4.15]{Bayraktar-Coman-Mariescu}. The inequality (\ref{ine-bichanmulogau}) is obtained in the estimate (4.16) in the proof of \cite[Lemma 4.15]{Bayraktar-Coman-Mariescu}.
\endproof

We note that the condition (\ref{ine-dktrenfzmu3}) is true 
if there is a number $C>0$ such that  $|f(z)| \le C(1+|z|)^{-3}$ for 
$z\in \C$ because
$$\int_{\{\log |z| >r\}} f d\Leb_\C 
\lesssim \int_{\{s >e^r\}}s^{-2} ds \lesssim e^{-r}.$$
Consequently, as a result of Theorem \ref{the-expectedzeros}, we obtain:
%===
\begin{theorem}[{\cite[Corollary 1.8]{Marinescu-Vu}}]
\label{the-expectedzeros-chuandensity} 
Let $(K,Q,\mu)$ be Bernstein-Markov. 
Assume that (H2) holds and there is a constant $C>0$ such that
$|f(z)| \le C(1+|z|)^{-3}$ for $z\in \C$.
Then we have 
\begin{align} \label{eq-corexpeczeros-chuandensity}
\E_k\big(k^{-1}[p_k=0]\big) = 
\ddc V_{K,Q}+ O\left(\frac{\log k}{k}\right).
\end{align}
%In particular if $C_k$ can be chosen to be equal to $O(k^s)$ for some constant $s\ in (0,1)$ independent of $k$, then the error current is equal to $O(k^{^-1+s})$. 
\end{theorem}

\subsection{Large deviation estimates}
We discuss now large deviation estimates which are stronger than Theorem \ref{the-expectedzeros-chuandensity}.

\begin{theorem}[A large deviaton type estimate I, {\cite[Theorem 1.7]{Marinescu-Vu}}]
\label{the-zeros2} Let the hypothesis be as in Theorem \ref{the-expectedzeros-chuandensity}.
Then for every constant $M \ge 1$ there exists a constant $C_M>0$ such that for every $k$,
\begin{align} \label{ine-largedevirationSZ2}
\mu_k\left\{(a_{k1}, \ldots, a_{kd_k}) \in \C^{d_k}: 
\dist_{-2}\big(k^{-1}[p_k=0], \ddc V_{K,Q}\big) \ge 
\frac{C_M\log k}{k} \right\} \le C_M k^{-M}.
\end{align}
 %where $\cali{P}_k$ denotes the joint-distribution of $\alpha_1,\ldots, \alpha_{d_k}$. 
\end{theorem}

We would also like to mention that in some cases, 
certain large deviation type estimates for random polynomials 
in dimension one were known; see \cite[Theorem 10]{Goetze-Jalowy} 
for polynomial error terms, and \cite[Theorem 1.1]{Dinh-random}, 
\cite[Theorem 3.10]{DV_random} for exponential error terms. 
To the best of our knowledge, there has been no quantitative 
generalization of the results in \cite{Shiffman-Zelditch} 
to higher dimension. It has been commented in the latter paper 
that their method does not seem to have a simple generalization 
to the case of higher dimension. 

%It was claimed in \cite[Theorem 1.7]{Marinescu-Vu} a somewhat cleaner statement but the proof given there only gives Theorem \ref{the-zeros}.

By \eqref{ine-interpolar}, one obtains similar estimates for 
$\dist_{-\beta}$ with $0< \beta \le 2$ as in Theorem \ref{the-zeros2}.  
The sharpness of the right side of \eqref{ine-largedevirationSZ2} is still unknown.

\begin{theorem}[A large deviaton type estimate II] \label{the-zeros-variant}
Let the hypothesis be as in Theorem \ref{the-expectedzeros-chuandensity}.  
Then there exists $A>0$ such that  for every number $\varepsilon >0$ 
there exists $C_\varepsilon>0$ satisfying
\begin{align} \label{ine-largedevirationSZvariant}
\mu_k\bigg \{(a_{k1}, \ldots, a_{kd_k}) \in \C^{d_k}: 
\dist_{-2}\big(k^{-1}[p_k=0], \ddc V_{K,Q}\big) \ge 
\varepsilon \bigg\} \le C_\varepsilon e^{-A \varepsilon k},
\end{align}
for every $k$. 
\end{theorem}

%where $O(\frac{\log k}{k})$ denotes a current $S_k$ 
%of order $0$ in $\C^n$ such that for every smooth form $\Phi$ 
%of degree $(2n-m)$ with compact support in $\C^n$ 
%such that $\|\Phi\|_{\mathcal{C}^2} \le 1$, we have 
%$$| \langle S_k, \Phi \rangle | \le C \frac{\log k}{k}\,,$$
%for some constant $C$ independent of $k, \Phi$. 
%===
The crucial ingredient in the proof of Theorem \ref{the-zeros2} and \ref{the-zeros-variant} is a polynomial growth of Bergman kernel functions, see Theorem \ref{the-convergencebergmankernel}.

Let $L$ be a complex algebraic curve in $\C^n$.
Since zero varieties of generic polynomials  intersect transversely $L$, almost surely the number of intersection 
points (without counting multiplicities) of the random hypersurface 
$\{p=0\}$ and $L$ is exactly $k\deg L$ by B\'ezout's theorem. Define 
$$\mu_{k,L}:= \frac{1}{k \deg L}\sum_{j=1}^{k \deg L} \delta_{z_j},$$
where $z_1, \ldots, z_{k \deg L}$ are zeros of $p$ on $L$. 
Let $[L]$ be the current of integration along $L$. 
Since $V_{K,Q}$ is bounded, the product 
$$\mu_L:= \frac{1}{\deg L}\ddc V_{K,Q} \wedge [L]$$
 is  a well-defined measure supported on $L$ 
 (it is simply $\ddc (V_{K,Q}|_{L})$ if $L$ is smooth). 

\begin{theorem}[{\cite[Theorem 1.9]{Marinescu-Vu}}]\label{the-zerosgiaoL} 
Let the hypothesis be as in Theorem \ref{the-expectedzeros-chuandensity}. 
Then for every constant $M \ge 1$ there exists a constant $C_M>0$ so that for every $k$,
$$\mu_k\left\{(a_{k1}, \ldots, a_{kd_k}) \in \C^{d_k}:
\dist_{-2}\big(\mu_{k,L}, \mu_L \big) \ge C_M\frac{\log k}{k}  \right\} 
\le C_M k^{-M}.$$ 
In particular, the measure $\mu_{k,L}$ converges weakly to 
$\mu_L$ almost surely as $k \to \infty$. 
\end{theorem}

%Now since the zero sets of $p_k$ on $L$ are discrete and equidistributed 
%with respect to $\mu_L$ as $k \to \infty$, one can ask as in \cite{Shiffman-Zelditch} 
%how they are correlated (if scaled appropriately). 
%Nevertheless such questions seem to be still out of reach in 
%the higher dimensional setting. Finally we note that one can 
%even consider $L$ to be a transcendental curve in $\C^n$. 
%In this case  generic polynomials $p$  still intersect $L$ 
%transversely asymptotically (see \cite{VietTuan}); 
%the issue of equidistribution is however more involved. 

Given that the zero sets of $p_k$ on $L$ are discrete and equidistributed 
with respect to $\mu_L$ as $k \to \infty$, one can pose the question of their 
correlation (if appropriately scaled) as in \cite{Shiffman-Zelditch}.
However, it appears that such inquiries remain beyond the reach 
in the higher dimensional setting. 
Finally, we note that we can even consider $L$ to be a transcendental curve in $\mathbb{C}^n$. 
In this case, generic polynomials $p$ still intersect $L$ transversely asymptotically 
(see \cite{VietTuan}). Nevertheless, the issue of equidistribution becomes more intricate.

In the last part of this subsection, we provide some explicit examples 
%(which are known in the literature) 
to which our main results apply. 

\begin{example} \label{exampleVK1} 
Let $K$ be the closed unit ball in $\C^n$. Then we have 
$V_K(z)= \log^+ \|z\|$ (where $\log^+ t:= \max\{\log t, 0\}$ for $t>0$); 
see \cite[Example 5.1.1]{Klimek}. 
Hence, we see that $V_K$ is Lipschitz but not continuously differentiable.   
\end{example} 

\begin{example}
Let $K=[-1,1]$ in $\C$. By \cite[Corollary 5.4.5]{Klimek}, we have 
$$V_K(z)=\log |z+ \sqrt{z^2-1}|$$ 
on $\C,$ where the square root is chosen such that 
$|z+ \sqrt{z^2-1}| \ge 1.$
In this case we have $V_K \in \mathcal{C}^{1/2}(X) \backslash \mathcal{C}^{1/2 
+\varepsilon}$ 
for any $\varepsilon \in (0,1/2).$ In higher dimension, similar observations also work 
for $K=[-1,1]^n \subset \C^n$ and $V_K$.   
\end{example} 
%===
\begin{example} 
Let $K$ be now the unit polydisk in $\C^n$. 
Thus by Example \ref{exampleVK1} and \cite[Theorem 5.1.8]{Klimek},
we have
$$V_K(z)= \max_{1 \le j \le n} \log^+|z_j|,$$ where
$z=(z_1,\ldots,z_n)$. Let $\mu$ be the restriction of the Lebesgue measure 
on $\C^n$ to $K$. For $J=(j_1,\ldots,j_n) \in \N^n$, put 
$J:= j_1+\cdots+j_n$, $z^J:= z^{j_1} \cdots z^{j_n}$. 
Observe that $(c_J z^J)_{J, |J| \le k}$ forms an orthonormal basis of 
$\mathcal{P}_k(K)$, where $c_J>0$ is a constant such that the norm
of $c_J z^J$ is equal to $1$. Let  
$$p_k:= \sum_{J, |J| \le k} \alpha_J c_J z^J,$$
where $\alpha_J$ are independent complex Gaussian random variables
of mean $0$ and variance $1$. Hence the hypothesis of
Theorem \ref{the-zeros2}
is fulfilled for $p_k$.
\end{example}

\subsection{Proof of quantitative equidistribution of random zeros}
\label{Subsec:Pqe}

In this subsection we prove Theorems \ref{the-zeros2} and its variants.  %We will first explain how to obtain Theorems \ref{the-BMBloomLevenberg} under an additional condition that $|f(z)| =O(|z|^{-3})$ as $z \to \infty$ in $\C$. 

Let $\Leb_{\C^m}$ be the Lebesgue measure on $\C^m$ for $m \ge 1$, and we denote by $\|\cdot\|$ the standard Euclidean norm on $\C^m$. Let $\omega_{FS,m}$ be the Fubini-Study form on $\P^m$, and let $\Omega_{FS, m}:= \omega_{FS,m}^m$ be the Fubini-Study volume form on $\P^m$. We always embed $\C^m$ in $\P^m$.

%\begin{lemma} (Borel-Cantelli lemma)
%Let $\{E_k\} \subset \mathcal{F}$ be a sequence of events on some probability space $(	\Omega,\mathcal{F}, \mathcal{P})$. If the sum of the probabilities of the $E_k$ is finite, i.e., 
%$$\sum_{k=1}^\infty \mathcal{P}(E_k) < \infty.$$
%Then the probability that infinitely many of them occur is 0, i.e.,
%$$\mathcal{P}(\cap_{k=1}^\infty \cup_{s \ge k} E_s)=0.$$
%\end{lemma}

%Let $f$ be a bounded Borel function on $\C$ such that there is a constant $C_0>0$ for which for every $r>1$ we have 
%$$\int_{|z| \ge r} f d \Leb_\C \le C_0/ r^2.$$ 
% Let $a_1,a_2,\ldots, a_{d_k} $ be complex-valued i.i.d  random variable whose distribution is $f \Leb_\C$. 

Let $p^{(d_k)}:= (s_1,\ldots, s_{d_k})$ be an orthonormal basis of $\mathcal{P}_k(\C^n)$.  Let $L$ be a complex algebraic subvariety of dimension $m \ge 1$ in $\C^n$. Note that the topological closure of $L$ in $\P^n$ is an algebraic subvariety in $\P^n$. %Thus the restriction of $\omega_{FS,n}^{m}$ to $L$ is a volume form on $L$. 
%Observe that 
%$$\omega_{FS,n}^m \le C \omega^m,$$
%where $\omega$ is the standard K\"ahler form on $\C^n$, and $C>0$ is a constant.    %Fix a compact $A$ of volume $\vol(A):=\int_A \omega_{FS,n}^m>0$ in $\C^n$. 
We start with a version of \cite[Lemma 2.4]{Bloom-Levenberg-random} with more or less the same proof. We use the Euclidean norm on $\C^{d_k}$.

\begin{lemma} \label{le-giongawn} Let $R \ge 1$ be a constant. 
Let $E_R$ be the set of $a=(a_1,\ldots, a_{d_k})\in \C^{d_k}$ so that $\|a\| \ge R$.
Then we have $$\mu_k(E_R) \le C  d_k^{2}R^{-1}$$
for some constant $C>0$ independent of $k$ and $R$.
\end{lemma}

We note that for $a \not \in E_R$, there holds
$$\log \bigg|\sum_{j=1}^{d_k} a_j s_j(z)\bigg| - 1/2\log \sum_{j=1}^{d_k} |s_j(z)|^2 \le \log R.$$

\proof  Observe that 
\begin{align*}
\mu_k(E_R) &\le \mu_k\bigg(a: |a_j| \ge d_k^{-1/2} R \, \text{ for some } 1\le j \le d_k\bigg) \\
&\le  d_k \int_{\{|a_1| \ge d_k^{-1/2}R\}} f \Leb_C \lesssim d_k^{3/2}R^{-1}.
\end{align*}  
\endproof

\begin{proof}[End of the proof of Theorems \ref{the-zeros2} and 
\ref{the-zerosgiaoL}] 
Recall that   $p_k= \sum_{j=1}^{d_k} a_{kj} s_j$. Put  $\psi_k:=
1/(2k) \log \widetilde{B}_k$.  Theorem \ref{the-convergencebergmankernel}  implies
$$\int_{L} | \psi_k -V_{K,Q}| \omega_{FS,n}^m \lesssim \frac{\log k}{k} \cdot$$
On the other hand,  observe that
$$\dist_{-2} \big(k^{-1}[p_k=0] \wedge [L], \ddc V_{K,Q}
\wedge [L]\big) \lesssim \int_{L} | k^{-1} \log |p_k| - V_{K,Q}|
\omega_{FS, n}^m.$$
It follows that 
\begin{align}\label{ine-distru3}
\dist_{-2} \big(k^{-1}[p_k=0] \wedge [L], \ddc V_{K,Q}
\wedge [L]\big) \lesssim \int_{L} | k^{-1} \log |p_k| - \psi_k|
\omega_{FS, n}^m+ O(k^{-1} \log k).
\end{align}
It remains to estimate $| k^{-1} \log |p_k| - \psi_k|$.
Let $M\ge 2$ be a constant and  $a=(a_{k1},\ldots, a_{kd_k}) \not \in E_R$ for $R:= d_k^{2M}$. By Lemma \ref{le-giongawn}, we obtain 
$$k^{-1} \log |p_k| - \psi_k \le 2M k^{-1}\log d_k, \quad \mu_k(E_R) \lesssim d_k^{-M}$$ 
Consequently,
\begin{align*}
|k^{-1} \log |p_k|- \psi_k| &= 2  \max\{k^{-1} \log |p_k|, 
\psi_k\}- \psi_k -k^{-1}\log |p_k| \\
&\le  \psi_k -k^{-1}\log |p_k|+ 2 M k^{-1} \log d_k.
\end{align*}
We note that 
$$- k(\psi_k -k^{-1}\log |p_k|) = \log \frac{\langle a, p^{(d_k)}\rangle}{ \|p^{(d_k)}\|}=: \xi(a).$$
For every constant $r>0$, put 
$$F_r:=\big\{ a \in \C^{d_k}: \xi(a) \le - r \big\}.$$
The inequality (\ref{ine-bichanmulogau}) applied to $u:= p^{(d_k)}$ tells us that 
$$\mu_k(F_r) \lesssim e^{-r}$$
Choose $r= M \log k$. We see that
$$\psi_k -k^{-1}\log |p_k| \le M k^{-1} \log k$$
for $a \not \in E_R \cup F_r$. Thus 
$$|k^{-1} \log |p_k|- \psi_k| \le 3M k^{-1} \log k$$
for $a \not \in E_R \cup F_r$.
 It follows that  
$$\int_{L} | k^{-1} \log |p_k| - \psi_k| \omega_{FS, n}^m 
\lesssim  k^{-1} \log k,$$
for $a \not \in E_R \cup F_r$. We already know that 
$$\mu_k(E_R \cup F_r) \lesssim k^{-M}.$$ 
This combined with (\ref{ine-distru3}) gives the desired assertion.
\end{proof}

\begin{proof}[End of the proof of Theorem \ref{the-zeros-variant}] The proof goes exactly as in the proof of Theorem \ref{the-zeros2}. Choosing $r:= \varepsilon k$ and $R:= e^{\epsilon k}$ and repeating arguments in the  proof of Theorem \ref{the-zeros2} gives
$$|k^{-1} \log |p_k|- \psi_k| \lesssim \epsilon $$
for $a \not \in E_R \cup F_r$, and 
$$\mu_k(E_R \cup F_r) \lesssim d_k^2 e^{-\epsilon k}.$$
\end{proof}

\subsection{Further discussions and questions}

We continue with some problems which 
have not been properly studied in higher dimensions.

\begin{problem} \label{pro-corre} Study the universal scaling limit 
of the correlations between zeros of random polynomials 
(e.g., when restricted to an algebraic curve $L$ as above).
\end{problem}

%As briefly mentioned above, Problem \ref{pro-corre} was investigated in \cite{Shiffman-Zelditch} 
%for dimension one. The method is based heavily on asymptotic of Bergman kernel functions. 
%We are not aware of any analogous work in the higher dimensional case except 
%\cite{BleherShiffmanZelditch} dealing with the situation where the non-pluripolar set $K$
%is equal to the whole manifold $\C^n$. 
%For a general $K$ as in the above hypothesis (H2), Problem \ref{pro-corre} 
%is still open for higher dimensions. 
%%In the next subsections we discuss the proof of Theorem \ref{the-zeros}.

As briefly mentioned above, Problem \ref{pro-corre} was investigated in 
\cite{Shiffman-Zelditch} for the dimension one. 
The method is heavily based on the asymptotic behavior of Bergman kernel functions. 
We are not aware of any analogous work in the higher dimensional case, 
except for \cite{BleherShiffmanZelditch}, which deals with the situation 
where the non-pluripolar set $K$ is equal to the entire manifold $\mathbb{C}^n$. 
For a general $K$ as in the above hypothesis (H2), Problem \ref{pro-corre} 
remains open for higher dimensions.

We come to the next problem. Theorem \ref{the-zeros2} only gives 
a large deviation type estimate. In the dimension one,  a much more precise estimate was proved in \cite{Zeitouni-Zelditch-largedeviation} (see also \cite{Butez_largedevia}) in the spirit of large deviation theory in probability. This implies the following estimate: for $n=1$ and every $\varepsilon>0$, there holds
\begin{equation}\label{ine-largeZZ}
\mu_k\Big\{(a_{k1}, \ldots, a_{kd_k}) \in \C^{d_k}:
\dist_{-2}\big(\mu_{k}, \mu \big) \ge \varepsilon  \Big\} 
\le C_\varepsilon \exp(-A_\varepsilon k^2),
\end{equation}
for some constants $C_\varepsilon>0,$  $A_\varepsilon>0$ and $k$ big enough. 
%The disadvantage in the above inequality compared to Theorem \ref{the-zeros} 
%is that $\delta$ must be independent of $k$. However the strong point is 
%the exponent rate $k^2$. We underline that the standard method using 
%the pluripotential theory has given so far only the exponent rate $k$
%in the estimate \eqref{ine-largeZZ}; see, e.g., \cite{Dinh-Ma-Marinescu}. 
%We are no aware of any higher dimensional generalization of results in 
%\cite{Zeitouni-Zelditch-largedeviation} at the moment of writing this survey. 
In contrast to Theorem \ref{the-zeros2}, the inequality \eqref{ine-largeZZ} 
has the disadvantage that $\varepsilon$ must be independent of $k$. 
However, its significant advantage lies in the exponent rate of $k^2$. 
It is important to note that the standard method employing pluripotential theory has, 
to the best of our knowledge, only provided the exponent rate $k$ in the estimate 
\eqref{ine-largeZZ}; for instance, see Theorem \ref{the-zeros-variant} above and \cite{Dinh-Ma-Marinescu}. 
Furthermore, we are not aware of any higher-dimensional generalizations of the 
results obtained in \cite{Zeitouni-Zelditch-largedeviation} at the time of writing this survey.

\begin{problem} (Large deviation for equidistribution of zeros)
Establish a higher dimensional generalization of the large deviation 
for equidistribution of zeros  in \cite{Zeitouni-Zelditch-largedeviation}.
\end{problem}

%An important consequence of the large deviation for equidistribution of zeros is an asymptotic for the probability of hole events, see \cite{Zeitouni-Zelditch-largedevia} and also \cite

To conclude this section, we want to discuss a recent exciting development on the distribution of zeros of random polynomials in dimension one. We are interested in the conditional equidistribution of zeros of random polynomials  in the event that zeros do not appear in a given subset of the ambient space. This problem was originally motivated by the problem of studying the distribution of physical gases when it is known that they avoid a certain set in space (\cite{Ghosh-Nishry-CPAM,Ghosh-Nishry}). Only some particular cases were treated for random point processes (in complex dimension one); see \cite{Ghosh-Nishry} for a survey. As far as we know, this is entirely open in higher dimensions. 

Similar questions for zeros of random sections of positive line bundles on compact Riemann surfaces were treated by  Dinh-Ghosh-Wu \cite{Dinh_Ghosh-Wu.holevent}, see also \cite{Wu-Xie.speedconvergence} for a refinement which gives a rate of converence for the conditional equidistribution. We refer to \cite{Drewitz-M-L3} for related results. It should be noted that the method in \cite{Dinh_Ghosh-Wu.holevent,Wu-Xie.speedconvergence}  relies again essentially in a variant of the notion of extremal plurisubharmonic functions, and an important part of their proofs is based on one-dimensional techniques. It is, thus, of great interest to   develop these methods to the higher dimensional setting.

\begin{problem} (Equidistribution of zeros in the hole event)
Extend the main results in  \cite{Dinh_Ghosh-Wu.holevent,Wu-Xie.speedconvergence} to higher dimensions.
\end{problem}

%\section{Zeros of holomorphic sections of 
%Hermitian vector bundles}
%
%\section{Toeplitz operators}
%
%\section{Strictly pseudocovex domains}

\bibliographystyle{siam}%{abbrv}

\bibliography{mmbook,biblio_family_MA,biblio_Viet_papers,bib-kahlerRicci-flow}
\end{document}